\documentclass{amsart}

\usepackage{epic}
\usepackage{eepic}
\usepackage{graphicx}
\usepackage{psfrag}

\newtheorem{thm}{Theorem}[section]
\newtheorem{ass}[thm]{Assumption}
\newtheorem{coro}[thm]{Corollary}
\newtheorem{lem}[thm]{Lemma}
\newtheorem{prop}[thm]{Proposition}

\theoremstyle{definition}

\newtheorem{eg}[thm]{Example}

\theoremstyle{remark}

\newtheorem{remk}[thm]{Remark}

\newcommand{\supp}{{\rm supp}}

\usepackage{graphicx}
\usepackage{psfrag}

\begin{document}

\title[Tracking eigenvalues to the frontier II]{
   \mbox{Tracking eigenvalues to the frontier  of moduli space II:} 
                         Limits for eigenvalue branches}

\author{Christopher M. Judge}       

\address{Indiana University, Bloomington, IN.}  

\email{cjudge@indiana.edu, 
          http://poincare.math.indiana.edu/\~{\/}cjudge}

\thanks{This work was supported in part by NSF 9972425.}



\date{\today}

\maketitle



\section{Introduction}
\label {SectionIntro}

Let $I \subset {\mathbb R}$ denote a compact interval symmetric about 
$0$ and let $c \cdot I$ denote the dilation of $I$ by a factor of $c$.  
Let $M$ and $N$, be compact oriented smooth manifolds of dimensions 
$d$ and $d+1$ respectively. We suppose that we have an embedded 
copy of $3 \cdot I \times M$ inside of $N$. (See Figure \ref{Degeneration}).
Let $N^0$ denote the complement of  the hypersurface $\{0\} \times M$. 

We consider families $\epsilon \rightarrow g(\epsilon)$
of Riemannian metric (tensors) on $N$ each of whose restriction
to $3 \cdot I \times M$ is a warped product of the following form
\begin{equation} \label{LocalMetric}
 g(\epsilon) |_{3 \cdot I \times M}~ 
 =~ \rho(\epsilon,t)^{2a}~ dt^2 
 + \rho(\epsilon,t)^{2b}~ h.
\end{equation}
Here $h$ is a fixed Riemannian metric on $M$,  
$\rho$ is smooth positive function that is positively homogeneous of degree $1$ on 
${\mathbb R}^2 \setminus \{\vec{0}\}$, and $a$ and $b$ are real numbers. 

\begin{figure}
  \centering
  \psfrag{a}{$a$}
  \psfrag{b}{$b$}
  \psfrag{+1}{$+1$}
  \psfrag{-1}{$-1$}
\psfrag{onezero}{hyperbolic \ 
                  degeneration}
  \psfrag{complete}{complete}
  \psfrag{incomplete}{incomplete}
  \psfrag{finitevolume}{finite volume}
  \psfrag{infinitevolume}{infinite volume}
    \includegraphics{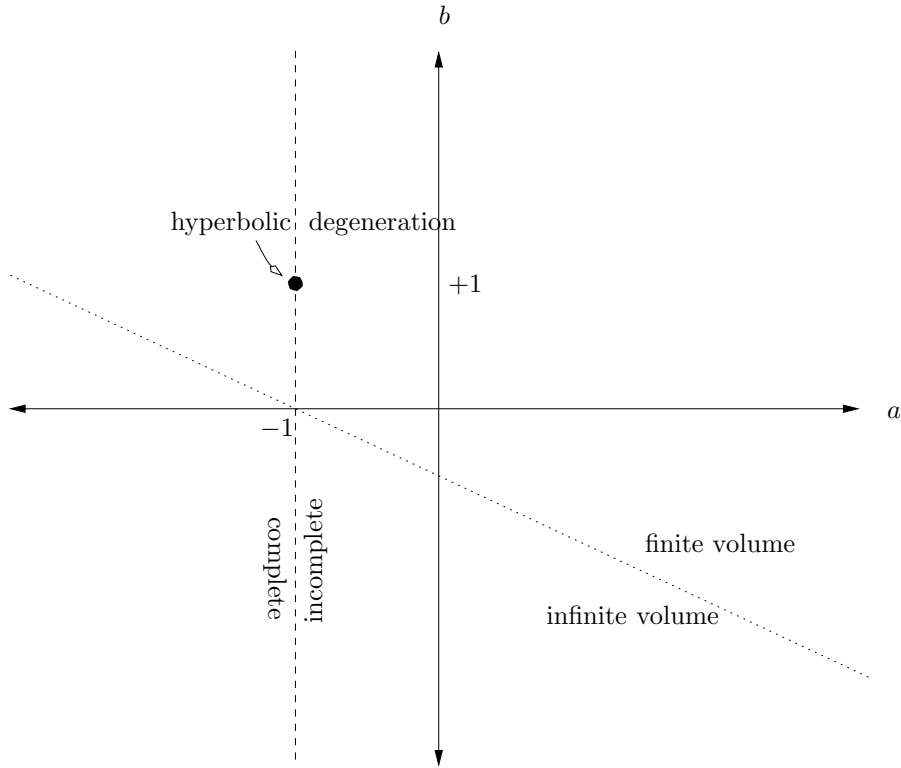}
  \caption{The $(a,b)$ parameter space}
  \label{Cases}
\end{figure}

The `limiting metric' $g(0)$ is singular 
along the hypersurface $\{0\} \times M$ provided
$(a,b) \neq \vec{0}$.  Indeed, since $\rho$ is
homogeneous 
\begin{equation*} 
  g(0) |_{3 \cdot I \times M}~ 
  =  ( c_{\pm} t)^{2a}~ dt^2 
  +~  ( c_{\pm} t)^{2b}~h.
\end{equation*} 
for some homogeneity constants $c_{\pm}$.
Melrose \cite{Mlr} has observed that the metric $g(0)$ 
is Riemannian complete if and only if $a \leq -1$, 
whereas $g(0)$ has finite volume if and only if $a+bd >-1$. 
(See Figure \ref{Cases}).

\medskip

\begin{eg}[Hyperbolic Degeneration]
\label{DegenExample}
Let $\gamma$ be a simple closed curve
in a compact oriented surface $N$ with $\chi(N) <0$. 
Let $g_{\epsilon}$ be a metric on $N$
of constant curvature $-1$ such that the unique geodesic 
homotopic to $\gamma$ has length 
$\epsilon< 2 \cosh^{-1}(2)$. By the collar lemma,
there exists an embedding 
$I \times \gamma \rightarrow N$ with
$I=[-\frac{1}{3}, \frac{1}{3}]$ such that
\begin{equation} \label{Classical}
  g_{\epsilon}|_{3 \cdot I\times \gamma}~ =~ 
  \frac{dt^2}{\epsilon^2 +t^2} 
  +  (\epsilon^2 +t^2)~ d x^2
\end{equation}
where $x$ is the usual coordinate on the circle ${\mathbb R} / {\mathbb Z} \cong \gamma$.
Note that the Riemannian surface $(I \times \gamma)^0,g_0)$ is 
a union of hyperbolic cusps. In this special case,
$\rho(\epsilon,t)= (\epsilon^2 + t^2)^{\frac{1}{2}}$, $a=-1$, and $b=1$.
\end{eg}

\medskip

Henceforth, we will assume that $\rho$ is strictly convex along
nonradial lines and that  $\partial_{\epsilon} \rho \geq 0$
for $\epsilon \geq 0$.  Moreover, we assume that the restriction of 
$g'(\epsilon)/g(\epsilon)$ to the unit tangent bundle of $K= N \setminus (I \times M)$ is 
bounded.\footnote{Here we view $g$ and $g'$ as quadratic forms on each tangent space; 
thus each may be regarded as function on the unit tangent bundle.} 
With these assumptions we have 

\medskip

\begin{thm}[Main Theorem] \label{MainFirst}
Let $b>0$, and either let $a < -1$ or let $a=-1$ and $a+bd=0$. 
Suppose that $\epsilon \rightarrow g(\epsilon)$ is a real-analytic family of 
Riemannian metrics on $N$ satisfying (\ref{LocalMetric}).  
Then each eigenvalue branch of the associated family of Laplacians, 
$\Delta_{g(\epsilon)}$, converges to a finite limit as $\epsilon$ tends  
to $0^+$.
\end{thm}

\medskip

\noindent
S. Wolpert \cite{Wlp92} proved Theorem \ref{MainFirst} in the special 
case of hyperbolic degeneration. He subsequently used this convergence
to produce evidence supporting the belief that Maass cusp forms `disappear' under
perturbation \cite{Wlp94} \cite{Snk95}. Note that although 
the present paper does not include a discussion of manifolds with cusps,
the methods described  here apply equally well to the eigenbranches 
of a `pseudo-Laplacian' associated to a manifold $N$ with finitely many cusps.

By combining results of this paper with those
of the prequel \cite{Jdg00}, we obtain
\begin{thm}
Let $a \leq -1$ and $b>0$.
Let $\psi(\epsilon)$ be an
eigenfunction branch whose
zeroth Fourier coefficient
(see \S \ref{SectionPreliminary})
vanishes identically for small 
$\epsilon$.
Then $\psi(\epsilon)$ converges 
to an 
$L^2(N^0,g(0))$-eigenfunction of 
$\Delta_{g(0)}$. 
\end{thm}

Since real-analytic eigenbranches can `cross',
their tracking is far subtler than the continuity 
of ordered eigenvalues.  
For example, consider the family of Laplacians, 
$\Delta_{\epsilon}$, associated to the flat tori
${\mathbb R}^2/ (\epsilon {\mathbb Z} \bigoplus 
  \epsilon^{-1} {\mathbb Z})$. In this case,
almost all of the real-analytic eigenvalue branches 
tend to infinity as $\epsilon$ tends to zero. Yet,
there are infinitely many branches that tend to zero.
Therefore if for each $\epsilon>0$, 
one were to label the eigenvalues in increasing order
(with multiplicities)
\begin{equation} \label{Label}
 0 < \lambda_1(\epsilon) 
  \leq \lambda_2(\epsilon) \leq \cdots \leq 
  \lambda_k(\epsilon) \leq  \cdots, 
\end{equation}
then each $\lambda_k(\epsilon)$ would tend to zero as $\epsilon$ 
tended to zero.\footnote{Note that although
$\lambda_k(\epsilon)$ is continuous, it is not real-analytic.}

Therefore, although (\ref{LocalMetric}) describes a relatively narrow
class of geometric degenerations, the conclusion of 
Theorem \ref{MainFirst} provides a great deal more information
concerning spectral behavior than the usual convergence 
results concerning ordered eigenvalues. (See, for example, 
the recent work of Cheeger and Colding \cite{ChgCld99}.)  
Indeed, the geometer's  standard tool for estimating the size of 
eigenvalues---the minimax principle---cannot be used 
to track real-analytic eigenvalue branches due to possible 
eigenbranch `crossings'.

Here, we rely instead on the variational principle 
$\dot{\lambda} = \int \psi \dot{\Delta}  \psi$. 
To illustrate our use of this principle, 
we prove in \S \ref{SectionGeneral} the following
general result:

\begin{figure}
  \centering
  \psfrag{t}{$t$}
  \psfrag{gamma}{$\gamma$}
  \psfrag{M}{$N$}
  \psfrag{U}{$I \times M$}
    \includegraphics{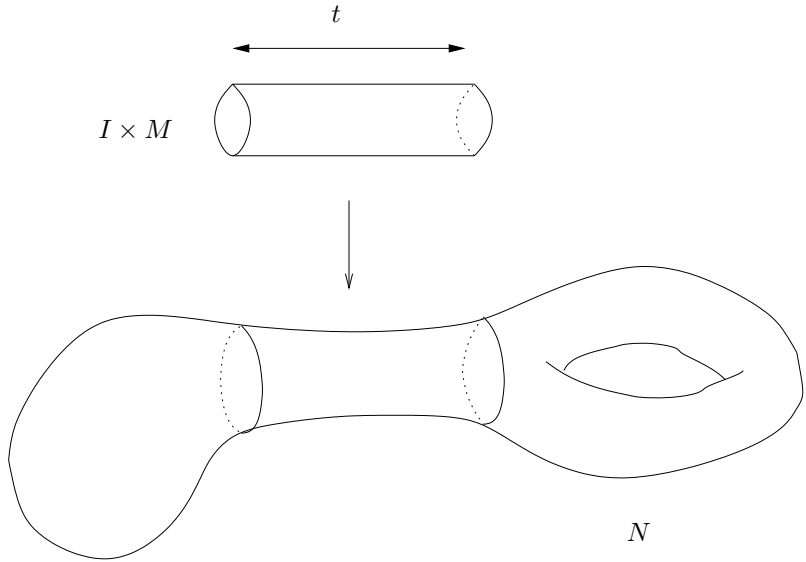}
  \caption{The set up}
  \label{Degeneration}
\end{figure}

\medskip

\begin{thm} \label{LogEstimatePropIntro}
Let $g(\epsilon)$ be a real-analytic 
family of metrics on a Riemannian manifold $N$. 
Then for each real-analytic
eigenbranch we have
\begin{equation} \label{EstGen}
  \left| 
  \frac{\lambda^{\prime}(\epsilon) }{
  \lambda(\epsilon)} \right|
  \leq  (\dim(N)+1) \left| \left|
  \frac{g'(\epsilon)}{g(\epsilon)}
  \right| \right|.
\end{equation}
\end{thm}

\medskip

For a family $g(\epsilon)$ satisfying (\ref{LocalMetric}), there
are vector fields supported in $I \times M$ such that 
the right hand side of (\ref{EstGen}) is unbounded 
as $\epsilon$ tends to zero.  Thus, in order to use 
the variational principle to prove Theorem \ref{MainFirst}, 
one must exhibit some control over the size of eigenfunctions 
in the bicollar $I \times M$.  For large eigenvalues,
controlling the size of eigenfunctions is notoriously difficult
\cite{Snk95} \cite{Zld00}. 
Indeed, for $b>0$, the central hypersurface $\{0\} \times M$ is 
totally geodesic, and hence the {\em correspondence principle}
of quantum physics leads one to `expect'---perhaps 
erroneously---that the mass of an eigenfunction with large eigenvalue
concentrates near $\{0\} \times M$.  The possibility of 
such `scarring' on $\{0\} \times M$ greatly contributes 
to the delicacy of the proof of Theorem \ref{MainFirst}.  
Fortunately, the ill-effects of possible `scarring' 
are ameliorated by the inequality $a \leq -1$,
that is, by the completeness of the limiting manifold.

We devote the remainder of this paper, with the exception of \S 
\ref{SectionBLessThan}, to proving Theorem  \ref{MainFirst}.
We now outline the contents and hence also the proof.  In \S
\ref{SectionGeneral} we illustrate our use of
the variational principle with a proof of 
Theorem \ref{LogEstimatePropIntro}.
In \S \ref{SectionPreliminary}, 
we establish some basic facts
concerning the warped product 
(\ref{LocalMetric}) including
an integration by parts formula
(Lemma \ref{IntByPartsLemma})
on which most of our analysis is based.

Underlying the proof of 
Theorem \ref{MainFirst} is a basic
fact: A nonnegative function 
$f \in C^1({\mathbb R}^+)$
has a finite limit as $\epsilon$
tends to $0^+$ provided 
the negative variation of $f$
over $]0, \epsilon]$ tends to zero as $\epsilon$ tends to zero.
Towards applying this to 
an eigenbranch $\lambda$, we 
derive  in \S \ref{SectionLowerBound}
lower bounds for 
the derivative $\lambda'$.
As an example of our approach, we 
use these lower bounds in
\S \ref{SectionBLessThan} to 
prove
\begin{thm} \label{B}
For $a < -1$ and $b\leq0$, each  
eigenvalue branch converges to a 
finite limit as $\epsilon$ tends  
to $0^+$. For $a = -1$ and 
$b <0$, each eigenvalue branch remains
bounded as $\epsilon$ tends to $0^+$.
\end{thm}
\noindent 
(Future work will include a 
more thorough investigation
of the cases in Theorem \ref{B}
as well as a study of the `adiabatic' 
case $(a,b)=(-1,0)$.)

Beginning with \S  \ref{SectionAPriori},
we restrict attention to the case of 
interest in the present work:
$a \leq -1$ and $b<0$.
We show in \S \ref{SectionAPriori} that 
$\epsilon^{2b} \cdot
   \lambda(\epsilon)$ converges to 
a finite limit 
(Theorem \ref{APrioriEstimate}). 
In \S \ref{SectionBootstrap}, we find 
that if $\epsilon^{2 k b} \cdot 
\lambda(\epsilon)$ remains bounded
for some $k<1$, then $\lambda(\epsilon)$
converges to a finite limit
(Theorem \ref{APrioriEstimate}). 
In \S \ref{SectionPsiBarThm}
boundedness for $k<1$
is verified provided 
$\mu^*= \lim \rho^{2b}(\epsilon,0)
 \cdot \lambda(\epsilon)$
is not a positive eigenvalue of 
the Laplacian $\Delta_h$
for $(M,h)$. Hence
in this case the eigenbranch has 
a finite limit (Theorem
\ref{ScarTheorem}).  
In \S \ref{SectionMainResult}
we assume that $\mu^*$ is a 
positive eigenvalue of  $\Delta_h$,
and obtain a contradiction in the
form of two conflicting estimates:
Lemmas \ref{LeftConflict} 
and \ref{RightConflict}.
Theorem \ref{MainFirst} follows.

We remark that the condition 
$\mu^* \in Spec(\Delta_h)
   \setminus \{0\}$---and hence
the threshold $k<1$---is intimately
tied to `scarring'. Indeed, 
one finds that the projection of 
$\psi$ onto the $\mu^*$-eigenspace 
is a scarring mode
in the sense of, for example, 
\S 7 of \cite{CdVPrs94}.
For the purpose of proving 
Lemma \ref{LeftConflict}
we need only know that the 
`width' of a scar is 
$O(\lambda^{-\frac{1}{4}})$
as $\lambda$ tends to infinity.
This result is given in Appendix
\ref{AppendixScar}.

The reader familiar with \S 3 in
\cite{Wlp92} will recognize the 
thread of the argument outlined above. 
Indeed, not only does 
the case of hyperbolic degeneration 
serve as motivation for the present work,
many of its basic features are 
representative of the general case.
On the other hand, at this
level of generality, we cannot 
avail ourselves of Teichm\"{u}ller 
theory nor the Poincar\'{e} series
estimate of \cite{Wlp92}.
Moreover, the peculiar features of 
the `overcomplete' case $a < -1$ 
do not appear in hyperbolic 
degeneration. These features add
complication to the arguments,
especially to those found in
\S \ref{SectionMainResult}.

I thank my ever-patient wife, Nacy, for her support. I also thank the referee 
for generous help wih the exposition.



\medskip

\section{Eigenvalue variation in the space of metrics}

\label{SectionGeneral}

\medskip

Let ${\mathcal M}(N) \subset 
             \otimes^2 TN^*$ be 
the space of all (smooth) Riemannian 
inner products on a compact manifold $N$.
To each $g \in  {\mathcal M}(N)$
we associate the Laplacian $\Delta_g$. 
This is a self-adjoint, unbounded 
operator on $L^2(N, dV_g)$ 
defined via the Friedrich's 
extension with respect to symmetric 
boundary conditions.

A fixed inner product 
$g^* \in  {\mathcal M}(N)$
induces a Banach norm on the space
of $2$-tensors $\otimes^2 TN^*$. 
A family of metric tensors 
$\epsilon \rightarrow g(\epsilon)$
is said to be real-analytic if 
it defines a real-analytic path 
in the Banach space $\otimes^2 TN^*$. 
Using the ratio of Riemannian measures
$dV_g/dV_{g^*}$,
one constructs a natural family 
of unitary operators that conjugates
$\epsilon \rightarrow
 \Delta_{g(\epsilon)}$ into 
a real-analytic family of compactly 
resolved operators 
that are self-adjoint with respect 
to the fixed sesquilinear form 
determined by $dV_{g*}$. 
It follows from analytic perturbation 
theory \cite{Kat} that there exists 
a countable collection of eigenfunction 
branches,  
$\{\epsilon \rightarrow 
   \psi_k(\epsilon) \} 
  \subset L^{2}(N,dV_g)
 \cap C^{\infty}(N)$,
such that for each fixed $\epsilon$,
the set $\{ \psi_k(\epsilon) \}$ 
is an orthonomal basis for  
$L^{2}(N, dV_g)$.

Given a continuous function 
$f: TN^* \setminus \{0\} 
      \rightarrow {\mathbb R}$ 
satisfying 
$f( c\cdot v)=f(v)$ for all 
$c \in {\mathbb R}$, 
let $||f||$ denote the supremum.  
An example of  such a function is 
$v \rightarrow h(v,v) / g(v,v)$ 
where $g \in {\mathcal M}(N)$ 
and $h \in \otimes^2 TN^*$ 
is an arbitrary 2-tensor.

\medskip

\begin{thm} \label{LogEstimateProp}
Let $g(\epsilon)$ be a real-analytic 
family of metrics on $N$. 
Then for each real-analytic
eigenbranch we have
\begin{equation} \label{LogEstimate}
  \left| 
  \frac{\lambda^{\prime}(\epsilon) }{
  \lambda(\epsilon)} \right|
  \leq  (\dim(N)+1) \left| \left|
  \frac{g'(\epsilon)}{g(\epsilon)}
  \right| \right|.
\end{equation}
\end{thm}
\medskip

\begin{proof}
We fix a background metric $g^*$
and write $dV^*$ for its volume form.
Define $\alpha \in C^{\infty}(N)$ by 
$dV_g(\epsilon)= \alpha(\epsilon) \cdot
dV^*$. 
Let $\psi(\epsilon)$ be an eigenfunction 
branch corresponding to 
$\lambda(\epsilon)$.
Supressing subscripts, we have
\begin{equation} \label{Minimax}
  \int  g(\nabla \psi, \nabla \psi) 
  \alpha \ dV^* 
  = \lambda \int \psi^2  \alpha  \ dV^*.
\end{equation}
Each object in (\ref{Minimax}) 
is real-analytic in $\epsilon$. 
By Taylor expanding, collecting first 
order terms, integrating by parts, 
and using the eigenequation, we find that
\begin{equation} \label{PreNormalize}
  \dot{\lambda} \int \psi^2 dV 
 = \int \dot{g}(\nabla \psi, 
 \nabla \psi) dV
+ 2 \int  g(\dot{\nabla} \psi, 
 \nabla \psi) dV
  +  \int   \frac{\dot{\alpha}}{\alpha} 
    \left( g(\nabla \psi, \nabla \psi) 
    - \lambda \psi^2 \right)  dV
\end{equation}
Here the symbol $\cdot$ denotes the first derivative with respect to $\epsilon$ 
evaluated at $\epsilon=0$.

By definition, we have 
$g(\nabla f, X)= X \cdot f$ for each 
fixed vector field
$X$ and function $f$ on $N$. 
Differentiating in $t$ yields
$\dot{g}(\nabla f, X) 
     + g(\dot{\nabla}f,X)=0$.
Using this identity, 
(\ref{PreNormalize}) reduces to 
\begin{equation} \label{WithGDot}
  \dot{\lambda} \int \psi^2 dV 
  = -\int \dot{g}(\nabla \psi, 
  \nabla \psi) dV
  + \int  \frac{\dot{\alpha}}{\alpha}
  \left( g(\nabla \psi, \nabla \psi) 
  - \lambda \psi^2 \right)  dV.
\end{equation}

From 
$|\dot{g}(\nabla \psi, \nabla \psi)
        /g( \nabla \psi, \nabla \psi)| 
      \leq ||\dot{g}/g||$, 
we have 
\begin{equation} \label{Comparison}
  \left| \int \dot{g}(\nabla \psi, 
  \nabla \psi) dV \right| 
  \leq ||\dot{g}/g|| \cdot \lambda
  \cdot \int  \psi^2 dV.
\end{equation}
By interpreting $\alpha(\epsilon)$ as 
the determinant of the matrix 
representation of $g(\epsilon)$ with 
respect to an orthonormal basis of $g^*$, 
one finds that the supremum of 
$|\dot{\alpha}|$ is bounded by 
$\dim(N) \cdot ||\dot{g}/ g||$. 
The claim follows by applying this 
bound and (\ref{Comparison}) 
to (\ref{WithGDot}).
\end{proof}

\medskip



\section{Preliminaries concerning 
              the warped product}
\label{SectionPreliminary}

\medskip

We record some basic facts concerning
the Laplacian, its eigenvalues, 
and eigenfunctions, on  $I \times M$ 
with the metric given in 
(\ref{LocalMetric}).
In the following $\Delta_h$, $\nabla_h$,
and $dV_h$, will denote respectively,
the Laplacian, gradient, and volume form, 
associated to the metric $h$ on
a fibre $\{t\} \times M$.

Recall that $d$ is the dimension of $M$.
For any $f \in C^{\infty}_0(I \times M)$  
\begin{equation} \label{Laplacian}
 \Delta_g f = - L(f)
 + \rho^{-2b} \Delta_h f 
\end{equation}
where 
\begin{equation}
 \label{L}
  L(f) = \rho^{-a-bd}~
             \partial_t~  \rho^{-a+bd}~
       \partial_t~ f.
\end{equation}
and
\begin{equation}
 \label{Gradient}
 \nabla_g f = 
 \rho^{-2a} \cdot \partial_t f
 + \rho^{-2b} \cdot \nabla_h f .
\end{equation}
The volume form restricted to $I \times M$ is 
\begin{equation}
  \label{Volume}
   dV_g = \rho^{a+bd}~dt~dV_h.
\end{equation}

\medskip

\begin{remk}
If no subscript appears, 
then the object is associated to $g$.
\end{remk}

\medskip

Given 
$f:I \times M \rightarrow {\mathbb R}$,
define
\begin{equation} 
  ||f||^2_M (t) 
    = \int_{\{t\} \times M}
       f^2(t,m)~ dV_h. 
\end{equation}

\medskip

\begin{prop} \label{BasicProp}
Let $\psi \in C^{2}(I \times M)$
satisfy $\Delta_g \psi =\lambda \psi$.
Then
\begin{equation} \label{LDiff}
  \frac{1}{2} L \left( ||\psi||^2_M \right)~
  = -\lambda \cdot ||\psi||^2_M~ +
     \int_{\{t\} \times M} 
     g(\nabla \psi, \nabla \psi)~ dV_h.
\end{equation}
\end{prop}

\medskip

\begin{proof}
Straightforward computation gives
\begin{equation*} 
  \frac{1}{2} L ( ||\psi||^2_M )~ 
  =~ \frac{1}{2} \int_{M} L (\psi^2) dV_h~ 
  =~   \int_{M} \psi L (\psi) dV_h 
  + \rho^{-2a} \int_{M}
 (\partial_t \psi)^2~ dV_h.
\end{equation*}
From  (\ref{Laplacian}) and $\Delta \psi = \lambda \psi$ we find that
\begin{equation*} \label{Energy}
  -\int_{M} \psi L (\psi)~ dV_h
  + \rho^{-2b} \int_{M}
 \psi \cdot \Delta_h \psi~ dV_h~
  =~ \lambda \int_M \psi^2~ dV_h. 
\end{equation*}
Integrating by parts over $M$
gives
\begin{equation}
\int_{M}\psi \cdot \Delta_h \psi~ dV_h
  = \int_{M}
    h(\nabla_h \psi, \nabla_h \psi)~ dV_h.
\end{equation}
Also note that from (\ref{Gradient})
we have 
\begin{equation} \label{GDecomposed}
  g(\nabla \psi, \nabla \psi) =
  \rho^{-2a} (\partial_t \psi)^2
  + \rho^{-2b}  
  h(\nabla_h \psi, \nabla_h \psi).
\end{equation}
The claim follows.
\end{proof}

\medskip

\begin{coro} \label{ConvexLemma}
Let $\psi \in C^{2}(I \times M)$
satisfy $\Delta_g \psi =\lambda \psi$.
Suppose that for each $t \in I$,
\begin{equation} \label{DirichletBound}
  \int_{\{t\} \times M} 
  ( \psi \cdot \Delta_h \psi)~  dV_h~
  \geq~ \mu \cdot \int_{\{t\} \times M} 
    \psi^2~ dV_h.
\end{equation}
Then
\begin{equation}
\frac{1}{2}  L \left( ||\psi||^2_M \right)~
  \geq~ \left(\mu \rho^{-2b}-\lambda \right)
 \cdot ||\psi||^2_M.
\end{equation}
\end{coro}

\medskip

\begin{proof}
Apply Proposition \ref{BasicProp}
and (\ref{GDecomposed}).
\end{proof}


%

\medskip

\begin{lem}[Integration by Parts Formula]
 \label{IntByPartsLemma}
Let $\sigma \in C^{\infty}({\mathbb R}^2 
        \setminus \{0\})$ 
be positive and positively homogeneous 
of degree $s$.  There exist constants 
$C, C'$, such that for any eigenpair
$(\psi, \lambda)$ on 
$(2 \cdot I) \times M$ we have 
\begin{equation*} 
\label{IntByPartsEstimate}
  \left|    \int_{I \times M} 
  \sigma \cdot \left( 
  g(\nabla \psi, \nabla \psi) 
  - \lambda  \psi^2 \right)~ dV \right|
  \leq  C \int_{I \times M} 
  \rho^{k} \psi^2 dV~
  +~ C'(\lambda+1)  
  \int_{(2I \setminus I) \times M} \psi^2~ dV
\end{equation*}
where $k=-2a-2+s$.
\end{lem}

\medskip

\begin{proof}
Let $J$ denote the dilated interval 
$\sqrt{2} \cdot I$, 
and let $0 \leq \chi \leq 1$
belong to $C_0^{\infty}(J)$ with 
$\chi \equiv 1$ on $I$.
By multiplying both sides of (\ref{LDiff})
by $\chi \sigma$ and integrating over 
$J$ one obtains
\begin{equation}
\int_{ J} \chi \sigma \cdot 
 L(||\psi||^2_M) \rho^{a+bd}dt
=   \int_{J \times M} 
 \chi \sigma \cdot \left( 
g(\nabla \psi, \nabla \psi) 
    - \lambda  \psi^2 \right)~ dV_g. 
\end{equation}
On the other hand, integration by parts gives
\begin{equation} \label{LSym}
\int_{J} \chi \sigma \cdot 
 L(||\psi||^2_M) \rho^{a+bd}dt
  = \int_{J} L(\chi \sigma) \cdot 
  ||\psi||^2_M \rho^{a+bd}dt.
\end{equation}
We have  $L(\chi \sigma)= \chi L(\sigma)
+ f$  where $f \equiv 0$  in a neighborhood
 of the origin. Note that $L$ adds $-2a-2$
to the homogeneity of any function.
Hence $\deg(L(\sigma))= -2a-2+s$,
and thus, since $\rho>0$,
we have  $|L(\sigma)| = O(\rho^k)$.
Therefore, by (\ref{LSym}) 
\begin{equation*}
\left| \int_{J} \chi \sigma \cdot 
 L(||\psi||^2_M) \rho^{a+bd}dt \right|~
\leq~ \int_{J \times M} \rho^{k} 
 \cdot \psi^2~ dV~ \leq~
 \int_{I \times M} \rho^{k} \psi^2 dV
 + C \int_{(2I\setminus I) \times M} \psi^2 dV
\end{equation*}
where $k=-2a-2+s$.

To complete the proof it suffices to show that 
\begin{equation} \label{Remainder}
  \left| \int_{J \times M} 
  (1-\chi)\cdot \sigma \cdot 
  g(\nabla \psi, \nabla \psi)~ dV \right|
  \leq  C'(\lambda +1) 
  \int_{(2 I \setminus I) \times M} \psi^2~ dV.
\end{equation}
Since $\sigma>0$ and $\sigma$ is 
bounded on the support of $1- \chi$, 
there exists $C'$ such that
\begin{eqnarray*}
   \int_{J \times M} 
   (1-\chi) \cdot \sigma \cdot 
   g(\nabla \psi, \nabla \psi)~ dV
   &\leq&  C' \int_{J \times M} 
     (1- \chi) 
    \cdot
   g(\nabla \psi, \nabla \psi)~ dV  
\end{eqnarray*}
Let $\eta(t)= (1-\chi(t)) \cdot 
                \chi(t/\sqrt{2})$. 
Integrating by parts gives
\begin{equation*}
   \int_{(2 \cdot I) \times M}  
    \eta 
   \cdot
   g(\nabla \psi, \nabla \psi)~ dV~
   =~   \lambda \int_{(2 \cdot I) \times M}
   \eta \cdot  \psi^2  dV
   + \int_{(2 \cdot I) \times M } \psi \cdot 
     g(\nabla \psi, \nabla \eta)~ dV.
\end{equation*}
Note that the support of $\nabla \eta$ belongs
to $(2 I \setminus I) \times M$.
Integration by parts in $t$ gives 
\begin{equation*}
   \int_{(2 I) \times M}  \psi \cdot 
   g(\nabla \psi, \nabla \eta)~ dV
   = \int_{(2I) \times M} 
    \psi^2 \cdot \partial_t^2  \eta~ dV.
\end{equation*}
Estimate (\ref{Remainder}) then 
follows from the fact that $\eta$ 
has support in $(2 \cdot I) \setminus I$  
and equals $1-\chi$ on $J$.
\end{proof}

\medskip

\begin{remk}  \label{HomoBound}
Suppose that $k \leq 0$. Then $\epsilon^{-k} \cdot \rho^k$ 
is well-defined, continuous, and homogeneous of degree 0. Hence it is bounded.
Thus, there exists a constant $C$ such that
\begin{equation}
  \rho^{k}(\epsilon, t)~ \leq~ C \cdot \epsilon^k
\end{equation}
for all 
$(\epsilon, t) \in {\mathbb R}^2 
        \setminus \{\vec{0} \} $.
\end{remk}

\medskip

\begin{remk} \label{ShiftMin}
Since $\rho$ is homogeneous and
strictly convex along nonradial lines, 
there exists $c \in {\mathbb R}$ 
such that for each $\epsilon \neq 0$,
the function 
$t \rightarrow \rho(\epsilon,t)$ 
has a unique maximum at $t= c \epsilon$.
To prove \ref{MainFirst}, without
loss of generality, we may assume that
$c=0$.  For  otherwise, 
based on the linear map $t \rightarrow
t -\epsilon c$, one may 
construct a real-analytic 
family of diffeomorphisms
$\phi_{\epsilon}: N \rightarrow N$
such that
\begin{equation}
 \phi_{\epsilon}^*g(\epsilon)|_{2I \times M}
 = \rho^{2a}(\epsilon, t-c\epsilon)~
 dt^2~
 +~ \rho^{2b}(\epsilon, t-c\epsilon)~ h. 
\end{equation}
Then one works with the 
positive, positively homogeneous
function $\rho(\epsilon, t-c\epsilon)$.
\end{remk}

\medskip

\begin{prop} \label{MaxedEpsilon}
For each $\epsilon \geq 0$, the 
maximum of the function 
\begin{equation} \label{MaxedOut}
 \sigma_{\epsilon}( t) 
 = \frac{\dot{\rho}}{\rho}~(\epsilon, t)~
\end{equation}
is $\epsilon^{-1}$. This maximum 
is uniquely achieved at $t=0$.
\end{prop}

\medskip

\begin{proof}
Note that by homogeneity and positivity, 
$\rho(\epsilon,0)= c \cdot \epsilon$ for some $c>0$, and hence 
$\sigma_{\epsilon}( 0) = \partial_{\epsilon}\log(\rho)(\epsilon,0) = \epsilon^{-1}$.
Therefore, the first claim will follow from
the second. 

By Remark \ref{ShiftMin}, we may assume that the
function $t \rightarrow  \rho(\epsilon, t)$ has a uniques minimum
at $t=0$, and thus $\rho^{-1}$ has a unique maximum there. 
To prove the claim, it will suffice to show the same for $\dot{\rho}$. 
In other words, it is enough to show that 
$\partial_{t}
           \partial_{\epsilon}
          \rho(\epsilon,t)$
is positive for $t<0$ and negative
for $t>0$.

Since $t \rightarrow \rho(\epsilon,t)$
is strictly convex, 
$t \rightarrow \partial_t \rho(\epsilon,t)$
is strictly increasing.
Let $0 \leq \epsilon_1 \leq \epsilon_2$.
The derivative
$\partial_t \rho$ is homogeneous of 
degree 0, and, therefore, for $t>0$
\begin{equation}
\partial_t \rho(\epsilon_2,t)
  =  \partial_t \rho \left(\epsilon_1, 
  \frac{\epsilon_1}{\epsilon_2} t \right)  
  \leq \partial_t \rho(\epsilon_1, t).
\end{equation}  
Hence $\epsilon \rightarrow 
   \partial_t \rho(\epsilon, t)$ is decreasing
for $\epsilon \geq 0$. That is,
$\partial_{\epsilon}
           \partial_t \rho(\epsilon,t)$
is negative for $t>0$ as desired.

An analogous argument shows that
$\partial_{\epsilon}
           \partial_t \rho(\epsilon,t)>0$
for $t<0$. The claim follows.
\end{proof}

\medskip




\bigskip

\section{Lower bounds for 
         degenerating families}

\label{SectionLowerBound}

\medskip

The purpose of this section is to derive a useful lower bound
for $\dot{\lambda}$.
Towards this end, we define the {\em zeroeth Fourier coefficient},
$\psi_0$, of a function $\psi$ on $I \times M$ by
\begin{equation} \label{Zeroeth}
    \psi_0(t)~ =~ \int_M \psi(t,m)~ dV_h(m)
\end{equation}
and the {\em complement}, $\widehat{\psi}$,  by
\begin{equation} \label{HatPsi}
   \widehat{\psi}(t,m)~ 
    =~ \psi(t,m) - \psi_0(t).
\end{equation}
Note that $\psi$ is a $\Delta_g$ eigenfunction with eigenvalue $\lambda$, 
if and only if both $\psi_0$ and $\widehat{\psi}$ are.

\medskip

In the sequel,  $K$ denotes the set complement $N \setminus (I \times M)$.

\medskip

\begin{thm} \label{MainLowerBound}
Let $a \leq 0$. There exist positive constants $C,C'$ 
such that for each eigenbranch $(\psi, \lambda)$  
\begin{eqnarray} \nonumber
  \dot{\lambda} \int_N \psi^2~ dV~ &\geq &
  -2 \max \{a,b\} \cdot \lambda \int_{I \times M} \frac{\dot\rho}{\rho}
  \cdot\widehat{\psi}^2~ dV~  \\  \label{MainLowerEstimate}
  && -~ C \int_{I \times M} \rho^{-2a-3} \cdot \psi^2~ dV~ 
  -~ C' (\lambda+1) \int_K \psi^2~ dV.
\end{eqnarray}
Moreover, if $a+bd=0$, then the integrand $\rho^{-2a-3} \cdot \psi^2$
can be replaced with  $\rho^{-2a-3} \cdot \widehat{\psi}^2$.

\end{thm}

\medskip

\begin{proof}
Our starting point is formula (\ref{WithGDot}):
\begin{equation}  \label{GDotSpecial}
    \dot{\lambda} \int_N \psi^2~ dV~
    =~ -\int_N \dot{g} (\nabla \psi, \nabla \psi)~ dV~
    +~  \int_N  \frac{\dot{\alpha}}{\alpha} \cdot
    \left( g(\nabla \psi, \nabla \psi) - \lambda \psi^2 \right)~  dV.
\end{equation}
Recall that by (global) hypothesis, the supremum 
of $|g'(\epsilon)/g(\epsilon)|$ over the unit tangent bundle of $K$ is finite.
By applying the argument that immediately follows (\ref{WithGDot})
to the restriction of  $g(\epsilon)$ to $TK$, we obtain 
\begin{equation} \label{GDotK} 
    \left|  -\int_K \dot{g}(\nabla \psi, \nabla \psi)dV
    + \int_K  \frac{\dot{\alpha}}{\alpha} 
    \left( g(\nabla \psi, \nabla \psi)
    - \lambda \psi^2 \right) dV \right|
    \leq C  \cdot \lambda  \int_K \psi^2 dV.
\end{equation}
for some positive constant $C$.
By (\ref{Volume}), the restriction of  $\dot{\alpha}/\alpha$ to $I \times M$
equals  $(a+bd)\dot{\rho}/\rho$. Thus, by combining  (\ref{GDotSpecial}) and
(\ref{GDotK}) we obtain
\begin{eqnarray}  \label{PreGDotEst} 
   \dot{\lambda} \int_N \psi^2~ dV  & \geq & \nonumber
   -\int_{I \times M}  \dot{g}(\nabla \psi, \nabla \psi)~dV~  \\
  & &  +~ (a+bd) \int_{I \times M} \frac{\dot{\rho}}{\rho}  \cdot
   \left( g(\nabla \psi, \nabla \psi) - \lambda \psi^2 \right)  dV~ \\
  & &  -~ C \cdot \lambda \int_{K} \psi^2~ dV.\nonumber
\end{eqnarray}

We claim that
\begin{equation} \label{HatReduction}
  \int_{I \times M} \dot{g}(\nabla \psi, \nabla \psi) dV~
  \leq~  \int_{I \times M} \dot{g}(\nabla \widehat\psi, \nabla \widehat\psi) dV.
\end{equation}
To see this, first note that from (\ref{LocalMetric}) we compute 
\begin{equation}  \label{GDotFormula}
  \dot{g}|_{I \times M}~  =~ \frac{\dot{\rho}}{ \rho} \cdot ( 2 a \rho^{2a} dt^2 + 2b \rho^{2b}h).
\end{equation}
The function $\psi_0$ is constant on each fibre $\{ t\} \times M$, and hence 
$h(\nabla \psi_0, \nabla \psi_0)=0$.  Therefore, since $a \leq 0$ and 
$\dot{\rho} \geq 0$,  we find that 
\begin{equation}  \label{Zeroth}
   \dot{g}(\nabla \psi_0,\nabla \psi_0)~
   =~ 2a \cdot \dot{\rho} \cdot \rho^{2a-1} \cdot 
   (\partial_t \psi_0)^2   \leq 0
\end{equation}
The operator $\partial$ preseves the decomposition
$\widehat{f} + f_0$.  In particular, $\int_M \partial \widehat \psi =0$ and
$\partial \psi_0$ is constant on each fibre. 
Therefore, $\int_M \partial \widehat{\psi} \cdot  \partial \psi_0~ dV_h =0$,
and it follows that
\begin{equation}  \label{CrossTerm}
  \int_{I \times M} 
  \dot{g}(\nabla \widehat{\psi}, \nabla \psi_0)~dV~
  =~ 0.
\end{equation} 
The claimed (\ref{HatReduction}) follows.

From (\ref{GDotFormula}) we also have that
\begin{equation} \label{XLevel}
 \dot{g}(X,X) \leq 2 \cdot \max\{a,b\}~ 
 \frac{\dot{\rho}}{\rho}~g(X,X)
\end{equation}
and hence combined with (\ref{HatReduction}) we have
\begin{equation} \label{Combo}
 \int_{I \times M} \dot{g}(\nabla \psi, \nabla \psi)~ dV~ 
  \leq~ 2 \cdot \max\{a,b\}~ \int_{I \times M}
 \frac{\dot{\rho}}{\rho}~g(\nabla \widehat\psi , \nabla \widehat \psi)~ dV.
\end{equation}
Substitution into (\ref{PreGDotEst}) then yields
\begin{eqnarray}  \label{PreIntEst}
   \dot{\lambda} \int_N \psi^2~ dV  & \geq & \nonumber
   - 2 \max\{a,b\} \int_{I \times M}  \frac{\dot{\rho}}{\rho} 
   \cdot g(\nabla \widehat \psi, \nabla \widehat \psi)~dV~  \\
  & &  +~ (a+bd) \int_{I \times M} \frac{\dot{\rho}}{\rho}  \cdot
   \left( g(\nabla \psi, \nabla \psi) - \lambda \psi^2 \right)  dV~ \\
  & &  -~ C \cdot \lambda \int_K \psi^2~ dV.\nonumber
\end{eqnarray}
Since $\dot{\rho}/\rho$ is homogeneous of degree $-1$,
Lemma \ref{IntByPartsLemma} applies to give
\begin{eqnarray} \label{AppliedIntByParts} \nonumber
  \left|  \int_{I \times M}  \frac{\dot{\rho}}{\rho} \cdot
  \left( g(\nabla \psi, \nabla \psi)
  - \lambda  \psi^2 \right)~dV \right|~
  &\leq~& C 
  \int_{I\times M} \rho^{-2a-3} \cdot \psi^2~ dV~  \\
  & &  +~  C'(\lambda + 1) \int_{(2I \setminus I)\times M}  \psi^2~ dV
\end{eqnarray}
as well as the analogous estimate with $\psi$ replaced by $\widehat \psi$. 
By combining these estimates with (\ref{PreIntEst}) and absorbing
constants, we obtain the claim.
\end{proof}

\medskip



\section{The case $b \leq 0$}

\label{SectionBLessThan}

\medskip

\begin{thm} \label{BLess}
Let $a < -1$ and $b \leq 0$. Then each eigenvalue branch 
$\lambda(\epsilon)$ converges to a finite limit as  $\epsilon$ tends to $0^+$.
\end{thm}

\medskip

\begin{proof}
We apply Theorem \ref{MainLowerBound}.
Since $\dot{\rho}$ and  $\rho$  are positive and $\max\{a, b\} \leq 0$,
the first term on the right hand side of (\ref{MainLowerEstimate}) is nonnegative.
Therefore, by Remark \ref{HomoBound} we have
\begin{equation}
 \dot{\lambda}~ \geq~ -C \epsilon^{-2a-3}~ -~ C' (\lambda +1)
\end{equation}
for some constants $C$ and $C'$.
Since $\lambda \geq 0$, division of both sides by $\lambda +1$ gives 
\begin{equation} \label{BNegEst}
 \frac{d}{d \epsilon} \log (\lambda +1)~
 \geq~  -C \epsilon^{-2a-3} -C'. 
\end{equation}
Since $a < -1$, the left hand side of  (\ref{BNegEst}) is integrable,
and, moreover, the negative variation of $\log (\lambda +1)$ over $]0, \epsilon[$
is $O(\epsilon^{-2a-2})$. Thus, since $-2a-2 >0$, the function $\log(\lambda +1)$ 
has a limit as $\epsilon$ tends to $0^+$.
Thus, the claim follows via exponentiation. 
\end{proof}

\medskip

\begin{prop} \label{BLess2}
Let $a \leq -1$ and $b < 0$.
Then each eigenvalue branch 
$\lambda(\epsilon)$ remains bounded 
as $\epsilon$ tends to $0^+$.
\end{prop}

\medskip

\begin{proof}
Let $\delta= -\max\{a, b\} >0$. Note that because $a \leq -1$, we have
$-2a-3 \geq -1$, and hence $\epsilon^{-2a-3} \leq  \epsilon^{-1}$
for $\epsilon$ small. Thus, by using Theorem \ref{MainLowerBound}
and Remark \ref{HomoBound}, we obtain
\begin{equation} \label{Ceiling}
  \frac{d}{d \epsilon} 
  \lambda~  \geq~ 
  (\delta \lambda - C) \cdot 
  \epsilon^{-1} - C'.
\end{equation}
If $\lambda > (C+1)/\delta$, then 
the right hand side of (\ref{Ceiling}) 
is positive for $\epsilon$ small.
The claim follows.
\end{proof}

\medskip


\section{An a priori estimate}

\label{SectionAPriori}


\medskip

\begin{ass} \label{Assump}
In the sequel we will assume that 
$b>0$ and either $a <-1$ or $a \leq -1$ and $a+bd=0$.
\end{ass}

\medskip


%

\begin{thm} \label{APrioriEstimate}
The quantity 
$\epsilon^{2b} \lambda(\epsilon)$
tends to a finite limit as $\epsilon$ 
tends to zero.
\end{thm}

\medskip

\begin{proof}
By Remark \ref{HomoBound} and Assumption \ref{Assump}, we have 
$\rho^{-2a-3} \leq C \cdot \epsilon^{-2a-3} \leq C \epsilon^{-1}$ for $\epsilon$ small.
Thus, it follows from Theorem \ref{MainLowerBound} and Proposition \ref{MaxedEpsilon}, 
that there exist positive constants $C$, $C'$ such that
\begin{equation}  \label{EpsilonBared2}
   \dot{\lambda}~
   \geq~  
-  2b \cdot \lambda \cdot \epsilon^{-1}~ -~ 
    C \cdot \epsilon^{-1}~  -~  C' \cdot\lambda~.
\end{equation}
Since $b>0$, we have upon 
letting $c= C/2b$
\begin{equation}
  \dot{\lambda}~ \geq~  
  -~ 2b \cdot (\lambda+c) \cdot  \epsilon^{-1}~ -~ C \cdot \lambda.
\end{equation}
Dividing by $\lambda +c$ gives
\begin{equation}  \label{EpsilonBared}
  \frac{d}{d \epsilon}
   \log{\left(\lambda +c \right)}~
  \geq~-  2b~ \cdot \epsilon^{-1}~ -~ C.
\end{equation}
Since  $\frac{d}{d \epsilon} \log(\epsilon^{2b}) = 2b \cdot \epsilon^{-1}$, 
we obtain 
\begin{equation} \label{Consequence}
\frac{d}{d \epsilon}
 \log{\left ( \epsilon^{2b} 
        (\lambda+c) \right)} \geq 
     -C.
\end{equation}
It follows that the negative variation of 
$f(\epsilon)=  \log{( \epsilon^{2b} (\lambda+c))}$
over the interval $]0, \epsilon[$ is $O(\epsilon)$.
It follows that $\lim_{\epsilon \rightarrow 0} f(\epsilon)$ 
is either finite, in which case 
$\lim_{\epsilon \rightarrow 0} \epsilon^{2b} \lambda(\epsilon)$
is finite, or  $\lim_{\epsilon \rightarrow 0} f(\epsilon)= - \infty$
in which case $\lim_{\epsilon \rightarrow 0} \epsilon^{2b} \lambda(\epsilon)=0$.
In either case the limit exists.
\end{proof}

\medskip



\section{A Bootstrap}

\label{SectionBootstrap}

Let $k_0$ denote the infimum of all $k$ such that the function 
$\epsilon^{2k b} \cdot \lambda(\epsilon)$ has a limit as $\epsilon$ tends to zero. 
By Theorem \ref{APrioriEstimate}, we have $k_0 \leq 1$.
The purpose of this section is to prove

\medskip


\begin{thm}   \label{Bootstrap}
If $k_0<1$, then $\lambda(\epsilon)$ tends to a finite
limit as $\epsilon$ tends to zero.
\end{thm}

\medskip

As a first step towards proving Theorem \ref{Bootstrap}, we 
have the following

\medskip

\begin{prop}
If there exist constants $C^*>0$ and $k<1$ such that 
\begin{equation} \label{InverseFractional}
 \int_{I \times M} \rho^{-1} 
 \widehat\psi^2~ dV~
 \leq~  C^* \epsilon^{-k} \int_{I \times M}  
 \widehat\psi^2 dV,
\end{equation}
then  $\lambda(\epsilon)$ tends to a finite
limit as $\epsilon$ tends to zero.
\end{prop}

\medskip

\begin{proof}
Since $\dot{\rho}$ is homogeneous of degree zero, it is bounded.  
By Assumption \ref{Assump}, we have $-2a-3 \leq -1$, and hence $\rho^{-2a-3} \leq  \rho^{-1}$.  
Therefore, via Theorem \ref{MainLowerBound} we find that
\begin{equation*}
   \dot{\lambda} \int_N \psi^2~ dV~
   \geq~
   - C \cdot (\lambda+1) \int_{I \times M}
  \rho^{-1} \cdot \widehat{\psi}^2~ dV~
  -~ C' (\lambda+1)
  \int_{K} \psi^2~ dV. 
\end{equation*}
Hence by (\ref{InverseFractional}) and Remark \ref{HomoBound}
\begin{equation}
  \dot{\lambda}~ 
  \geq~ - C'' \cdot (\lambda+1)  \cdot \epsilon^{-k}
\end{equation}
for some positive $C''$ and $\epsilon$ small.
Dividing by $\lambda +1$ gives
\begin{equation}
  \frac{d}{d \epsilon} \log \left(
   \lambda + 1 \right)~ \geq~ - C'' \cdot \epsilon^{-k}. 
\end{equation}
Since $k<1$, the right hand side is integrable,
and, in particular, the negative variation  of 
$f(\epsilon)=\log(\lambda + c)$ is $O(\epsilon^{1-k})$. 
Therefore $\lim_{\epsilon \rightarrow 0} f(\epsilon)$ exists,
and it follows that $\lambda(\epsilon)$ has a limit.
\end{proof}

\medskip

To verify (\ref{InverseFractional})---and thus prove Theorem \ref{Bootstrap}---we 
split the domain of integration of the integral on the left hand
side according to whether Corollary \ref{ConvexLemma}
implies the convexity of $||\psi||^2_M$ or not.
To be precise, let $\mu_1$ denote the smallest non-zero 
eigenvalue of $\Delta_h$. Define $A(\epsilon)$ to be set of $t$ such that 
\begin{equation} \label{AInequality}
 \lambda(\epsilon) \cdot \rho^{2b}(\epsilon, t)~
   \leq~  \frac{\mu_1}{2}.
\end{equation}
The key idea in what follows is that (\ref{AInequality}) and 
Corollary \ref{ConvexLemma} imply that the function $t \rightarrow ||\psi||^2_M(t)$
is convex enough to tame the singular behavior
of $\frac{\dot{\rho}}{\rho}$ near $(0,0)$. This heuristic
will be made precise in Lemma \ref{BootstrapLemma}.

In the following $B(\epsilon)$ will denote the set complement 
$I \setminus \frac{1}{2} A(\epsilon)$.

\medskip

\begin{proof}[Proof of Theorem \ref{Bootstrap}]
We claim that if $k_1>k_0$, then for sufficiently small $\epsilon$
\begin{equation} \label{BEstimate}
 \int_{B(\epsilon) \times M}
 \rho^{-1} |\widehat\psi|^2~ dV~ \leq~ 
 \epsilon^{-k_1}  \int_{I \times M} 
 |\widehat\psi|^2 dV_g.
\end{equation}
Indeed, from (\ref{AInequality}), we have
$t \in B(\epsilon)$ if and only if   
\begin{equation}
  \lambda(\epsilon) \cdot \rho^{2b}(\epsilon,2t)~ \geq~ \frac{\mu_1}{2}.
\end{equation}
If $k_1>k_0$, then for sufficiently small $\epsilon$, 
we have $\epsilon^{2bk_1} \cdot \lambda(\epsilon) \leq \frac{\mu_1}{2}$,
and hence for $t \in B(\epsilon)$ 
\begin{equation}
  \rho^{2b}(\epsilon,2t)~ \geq~  \epsilon^{2bk_1}.
\end{equation}
Therefore $\rho^{-1}( 2 \epsilon, 2t) \leq (2 \epsilon)^{-k_1}$
for $t \in B(\epsilon)$ and sufficiently small $\epsilon$. 
The claimed (\ref{BEstimate}) then follows from homogeneity and integration.

By applying Lemma \ref{BootstrapLemma} with $\mu=\mu_1$, $\delta= \frac{\mu}{2}$,
and  $\psi= \widehat{\psi}$, 
we obtain the complementary estimate.
Indeed, since $M$ is compact, the function 
$\widehat{\psi}$---defined by (\ref{HatPsi})---is orthogonal 
to the zero eigenspace of $\Delta_h$. Hence for each $t \in I$
\begin{equation}
  \int_{\{t\} \times M}
   \widehat\psi \cdot \Delta_h  \widehat\psi~ dV_h~
 \geq~ \mu_1 \int_{\{t\} \times M}  
  |\widehat\psi|^2 dV_h.
\end{equation}
Thus, Lemma \ref{BootstrapLemma} applies to give
\begin{equation} \label{AEstimate}
\int_{\frac{1}{2}A(\epsilon) \times M}
\rho^{-1} |\widehat\psi|^2~ dV_g~ \leq~ 
  C \cdot \epsilon^{-(1-2b)}  \int_{I \times M} 
 |\widehat\psi|^2 dV_g.
\end{equation}
Combining (\ref{BEstimate}) and (\ref{AEstimate})
gives us the desired (\ref{InverseFractional}) for all
$k > \max\{1-2b, k_0\}$.
\end{proof}

\medskip



\begin{remk} \label{RhoLess}
Without loss of generality, we may assume that $\rho <1$ on $I \times I$.
\end{remk}

\medskip

\begin{lem} \label{BootstrapLemma}
Let $\mu>\delta>0$ and let $(\psi, \lambda)$ be a Laplace eigenbranch on $I \times M$.
Let $A(\epsilon)$ denote the set of $t$ that satsify
\begin{equation} \label{ConvexInequality}
 \lambda(\epsilon) \cdot \rho^{2b}(\epsilon, t)~
   \leq~  \mu - \delta.
\end{equation}
Suppose that 
\begin{equation} \label{Dies}
\lim_{\epsilon \rightarrow 0}  \lambda(\epsilon)  \cdot \rho^{2b}(\epsilon,0)~ <~ \mu-\delta
\end{equation}
and that for each $t$
\begin{equation} \label{OrthogonalConstantsLowerBound}
 \int_{\{t\}\times  M} \psi \cdot \Delta_h \psi~ dV_h~
  \geq~ \mu  \int_{\{t\}\times  M}  \psi^2~ dV_h.
\end{equation}
Then there exist $C>0$ such that for small $\epsilon>0$
\begin{equation} \label{The52Estimate}
  \int_{\frac{1}{2} A(\epsilon) \times M} 
  \rho^{-1} \cdot \psi^2~ dV_g~
  \leq~  C  \cdot \epsilon^{2b-1} 
  \int_{A(\epsilon) \times M} 
        \psi^2~ dV_g.
\end{equation}
\end{lem}

\medskip

\begin{proof}
By Remark \ref{RhoLess} we have $\delta \rho^{2b} \leq \delta$
and hence it follows from (\ref{ConvexInequality}) that
$\delta \leq  \mu \cdot \rho^{-2b} - \lambda$. Therefore
\begin{eqnarray}
\nonumber
 \delta   \int_{\frac{1}{2}A(\epsilon) \times M}
  \rho^{-1} \psi^2~ dV  &=&
 \delta   \int_{\frac{1}{2} A(\epsilon)} 
   \rho^{-1}
  ||\psi||^2_M \cdot  \rho^{a+bd}~ dt   \\
  &\leq& \int_{\frac{1}{2} A(\epsilon)} 
  (\mu \rho^{-2b}-\lambda) \cdot 
  ||\psi||^2_M  \cdot 
   \rho^{2b-1} \cdot
  \rho^{a+bd} dt.  \nonumber  \\
  &\leq& 
   \frac{1}{2} \int_{\frac{1}{2} A(\epsilon)}
    L(||\psi||^2_M) \cdot 
  \rho^{2b-1} \cdot \rho^{a+bd} dt.
  \label{PreParts} 
\end{eqnarray}
Here the last inequality
follows from (\ref{OrthogonalConstantsLowerBound}) and Corollary \ref{ConvexLemma}.

We wish to apply integration by parts to the last integral in 
(\ref{PreParts}). Towards this end, let $\chi$ be a smooth function
supported in  $[-1,1]$ with $\chi \equiv 1$ on $[-1/2,1/2]$ and
$\max |\chi|= 1$.  Since $\rho$ is positive, convex, and homogeneous of 
degree 1, the set $A(\epsilon)$ is a closed interval 
$[t_-(\epsilon), t_+(\epsilon)]$ that contains $0$.
For each $\epsilon>0$ define
\begin{equation}
    \bar{\chi}(\epsilon,t) = \left\{ \begin{array}{ll}
    \chi \left(\frac{t}{t_+(\epsilon)} \right) & \rm{for} \ t \geq 0 \\
    \chi \left(\frac{t}{t_-(\epsilon)} \right) & \rm{for} \ t \leq 0. \end{array} \right.
\end{equation}

Integration by parts shows that the operator $L$
is symmetric on $L^2(I, \rho^{a+bd}~dt)$ with Dirichlet
boundary conditions. Thus, 
\begin{equation} \label{PartsForL}
 \int_{\frac{1}{2} A(\epsilon)} \bar{\chi} \cdot
  L(||\psi||^2_M) \cdot 
  \rho^{2b-1} \cdot \rho^{a+bd} dt
=  \int_{I}   ||\psi||_M^2  \cdot
   L(\bar\chi \cdot \rho^{2b-1}) \cdot \rho^{a+bd}~ dt.
\end{equation}
By (\ref{OrthogonalConstantsLowerBound}) and
Corollary  \ref{ConvexLemma}, we have  $L(||\psi||^2_M(t)) >0$
for $t \in A(\epsilon)$. Thus, since $A(\epsilon)= \supp(\bar{\chi})$,
estimate (\ref{PreParts}) and (\ref{PartsForL}) combine to give
\begin{equation}
2 \delta \int_{\frac{1}{2}A(\epsilon) \times M}
  \rho^{-1} \psi^2~ dV~
\leq~  \int_{I \times M}     
   L(\bar{\chi} \cdot \rho^{2b-1})\cdot  \psi^2~ dV.
\end{equation}
Therefore, to verify (\ref{The52Estimate}), it will suffice to show that
\begin{equation} \label{LEstimate}
  \left| L(\bar\chi \cdot \rho^{2b-1}) \right| 
  = O( \epsilon^{2b-1}).
\end{equation}

By homogeneity, the supremum of $\rho^c(\epsilon, t)$ 
over $I$ is $O(\epsilon^c)$ for any constant $c$. We compute 
\begin{equation} \label{LCompute}
  L(\bar\chi \rho^{2b-1})~
  =~ \bar\chi \cdot L(\rho^{2b-1})~
  +~ \beta \cdot \rho^{-2a +2b -2} 
  \partial \rho~ \partial \bar\chi~
+~ \rho^{-2a+2b-1} \partial^2 \bar \chi
\end{equation}
where $\beta= (-a +bd +2c)$.
The operator $L$ adds $-2a-2$ 
to the degree of a homogeneous function,
and hence 
\begin{equation} \label{Desire1}
|L(\rho^{2b-1})| \leq  \rho^{2a-3+2b}
 =O(\epsilon^{2a-3+2b})
\end{equation}
Since $a \geq -1$, we have $2a-3 \geq -1$, 
and hence $|L(\rho^{2b-1})|= O(\epsilon^{2b-1})$.

To estimate the remaining two terms in (\ref{LCompute}), we need to estimate
$|\partial \bar \chi|$ and $|\partial^2 \bar \chi|$.
To this end, consider  $r(\epsilon)= \min \{|t_{\pm}(\epsilon)|\}$, 
the inner radius of $A(\epsilon)$. By Lemma \ref{DiameterBound}, there exists  $\eta>0$
such that $|r(\epsilon)| \geq \eta \cdot \epsilon$ for all $\epsilon$ small. 
Therefore  $|\partial \bar \chi| \leq r(\epsilon)^{-1} \cdot |\chi'|= O(\epsilon^{-1})$ 
and, similarly, $|\partial^2 \bar \chi| = O(\epsilon^{-2})$.

The function $\rho^{-2a +2b -2}$ appearing in (\ref{LCompute}) is 
homogeneous of degree $-2a-2 +2b \geq 2b$. Hence since
$\partial \rho$ is homogeneous of degree  zero,
\begin{equation} \label{2Desire}
 | \rho^{-2a +2b -2}  \cdot
  \partial \rho \cdot \partial \bar \chi|
  = O(\epsilon^{2b-1}).
\end{equation}
A similar argument shows that
\begin{equation}  \label{3Desire}
 | \rho^{-2a +2b -1} \cdot
  \partial \rho \cdot \partial^2  \bar \chi|
  = O(\epsilon^{2b-1})
\end{equation}
The desired estimate (\ref{LEstimate}) then follows from (\ref{Desire1}),
(\ref{2Desire}), and (\ref{3Desire}).
\end{proof}

\medskip

\begin{lem} \label{DiameterBound}
Let $(\psi, \lambda)$ be a Laplace eigenbranch on $I \times M$.
Let $\mu>\delta>0$ and let  $A(\epsilon)=[t_-(\epsilon), t_{+}(\epsilon)]$ be defined 
as in (\ref{ConvexInequality}).  If (\ref{Dies})
holds, then there exists $\eta>0$ such that
\begin{equation}
     \left|t_\pm(\epsilon) \right| \geq \eta \cdot \epsilon.
\end{equation}
for small $\epsilon$.
\end{lem}

\medskip

\begin{proof}
By definition, $t_{\pm}$ satisfies
\begin{equation*}
    \lambda(\epsilon) \cdot \rho^{2b}( \epsilon,  t_{\pm}(\epsilon))~ =~ \mu- \delta.
\end{equation*}
Thus, by homogeneity
\begin{equation} \label{Implicitly}
    \lambda(\epsilon)  \cdot \epsilon^{2b} \cdot \rho^{2b}( 1,  \epsilon^{-1} \cdot 
           t_{\pm}(\epsilon))~ =~ \mu-\delta.
\end{equation}
By Theorem \ref{APrioriEstimate},  $\lambda(\epsilon)  \cdot \epsilon^{2b}$ 
tends to $c \leq 0$ as $\epsilon$ tends to zero. If $c=0$, then since 
the right hand side of (\ref{Implicitly}) is positive,  
$\epsilon^{-1} \cdot t_{\pm}(\epsilon)$ must tend to infinity
as $\epsilon$ tends to zero.  

By homogeneity, (\ref{Dies}) implies that  $c \cdot \rho^{2b}(1,0)< \mu-\delta$.
Thus if $c>0$, then by (\ref{Implicitly}) we have
\begin{equation*} 
   \rho^{2b}( 1,  \epsilon^{-1} \cdot t_{\pm}(\epsilon))~ 
   =~ \frac{\mu-\delta}{c}~ >~ \rho^{2b}(1,0)~ 
\end{equation*}
By Remark \ref{ShiftMin} and the condition $b>0$,  the function 
$t \rightarrow \rho^{2b}(1,t)$ assumes a unique minimum at $t=0$.
Therefore, $\epsilon^{-1} \cdot t_{\pm}(\epsilon)$ is strictly bounded away
from zero. The claim follows.
\end{proof}

\medskip



\section{A second bootstrap}

\label{SectionPsiBarThm}

\medskip

By homogeneity $\rho^{2b}(\epsilon, 0)= \rho^{2b}(1, 0) \cdot \epsilon^{2b}$,
and hence, by Theorem  \ref{APrioriEstimate}, the limit
\begin{equation}
\label{mDefn}
\mu^*~=~ \lim_{ \epsilon \rightarrow 0}~ 
\rho^{2b}(\epsilon, 0) \cdot 
             \lambda(\epsilon)
\end{equation}
exists. The purpose of this section is to prove 

\medskip

\begin{thm} \label{ScarTheorem}
If $\mu^*$ does not equal a positive eigenvalue of $\Delta_h$, then 
$\lambda(\epsilon)$ tends to a finite limit as $\epsilon$ tends to zero.
\end{thm}

\medskip

To prove Theorem \ref{ScarTheorem}, we will use the $\Delta_h$ spectral
decomposition of $\psi$. To be precise, the orthogonal projection,  
$E_{\mu}: L^2(M,dV_h) \rightarrow  L^2(M,dV_h)$, 
onto the $\mu$-eigenspace of $\Delta_{h}$ extends fibre by fibre to an operator 
$\bar{E}_{\mu}: L^2(I\times M ,dV_g) \rightarrow L^2(I\times M ,dV_g)$.  Set
\begin{equation}  \label{PsiStar}
   \psi^* =  \bar{E}_{\mu*} (\psi) 
\end{equation}
\begin{eqnarray*} 
  \psi_+ &=& \sum_{\mu>\mu^*} 
  \bar{E}_{\mu} (\psi) \\
  \psi_-&=&  
  \sum_{0<\mu<\mu^*} \bar{E}_{\mu} (\psi).
  \nonumber
\end{eqnarray*}
Note that $\psi^*$, $\psi_+$ and 
$\psi_-$ are all eigenfunctions 
of $\Delta_g$ with eigenvalue $\lambda$. 
The $0$-eigenspace of $\Delta_h$ 
consists of the constant functions. 
Thus if $\mu^*$ is not a positive eigenvalue, 
then 
\begin{equation} \label{NoStar}
\psi= \psi_0 + \psi_{-} + \psi_{+}.
\end{equation}
where $\psi_0$ is defined in (\ref{Zeroeth}).

\medskip

\begin{proof}[Proof of Theorem  \ref{ScarTheorem}]
Since $a<0$, by (\ref{GDotFormula}) we have
\begin{equation*}
  \dot{g}(\nabla \widehat\psi, \nabla \widehat\psi)~ 
  \leq~ 2b \cdot \frac{\dot{\rho}}{\rho} \cdot \rho^{2b} \cdot 
   h(\nabla_g \widehat\psi, \nabla_g \widehat\psi)
\end{equation*}
Note that from (\ref{Gradient}) we obtain
$h(\nabla_g\widehat \psi, \nabla_g \widehat\psi)
   = \rho^{-4b} h(\nabla_h \widehat\psi, \widehat\nabla_h \psi)$.
Hence,
\begin{equation} \label{bPosBound}
 \int_{I \times M}  \dot{g}(\nabla \widehat\psi, \nabla \widehat\psi)~ 
  \leq~ 2b \int_I \frac{\dot{\rho}}{\rho} \cdot \rho^{-2b} 
    \left(\int_{\{t\} \times M} h(\nabla_h \widehat\psi, \nabla_h \widehat\psi)~ dV_h \right)~ dt.
\end{equation}

By hypothesis we have $\psi^*=0$ and hence $\widehat{\psi}= \psi_+ + \psi_-$.  
Applying Parseval's principle for $\Delta_h$ acting on $L^2(M, dV_h)$, we obtain
\begin{equation} \label{hParseval}
  \int_{\{t\} \times M} h(\nabla_h \widehat\psi, \nabla_h \widehat\psi)~ dV_h~
  =~  \int_{\{t\} \times M} h(\nabla_h \psi_-, \nabla_h \psi_-)~ dV_h~
  +~  \int_{\{t\} \times M} h(\nabla_h \psi_+, \nabla_h \psi_+)~ dV_h
\end{equation}
as well as
\begin{equation} \label{PlainParseval}
  \int_{\{t\} \times M} \psi^2~ dV_h~
  =~  \int_{\{t\} \times M} (\psi_-)^2~ dV_h~
  +~  \int_{\{t\} \times M} (\psi_+)^2~ dV_h.
\end{equation}
Thus, it follows from Lemmas \ref{MinusBound} and \ref{PlusBound} that 
there exists $k<1$ such that
\begin{eqnarray} \label{AppliedhParseval} \nonumber
  \int_{I \times M} \frac{\dot{\rho}}{\rho} \cdot \rho^{-2b} 
    \cdot  h(\nabla_h\widehat \psi, \nabla_h \widehat\psi)~ dV~
  &\leq& \left( k \cdot \lambda+ C \epsilon^{-2a-2} \right)
  \cdot \epsilon^{-1} \int_{I \times M}\widehat \psi^2~ dV~ \\
 & & +~    C' \cdot (\lambda+1) \int_K \psi^2~ dV.
\end{eqnarray}
Since $-2a-2 \geq 0$, combining (\ref{AppliedhParseval})
with (\ref{bPosBound})  and (\ref{HatReduction}) yields $k<1$ such that
\begin{equation} \label{PsiBarEstimate}
 \int_{I \times M} \dot{g}(\nabla \psi, \nabla \psi)~ dV~ \leq~
 2b \cdot  k \cdot (\lambda+c) \cdot  \epsilon^{-1}   
 \int_{I \times M} \widehat\psi^2~ dV~ +~ C'\cdot (\lambda+1) \int_K \psi^2 dV
\end{equation}
for some $c>0$.

By substituting (\ref{PsiBarEstimate}) into (\ref{PreGDotEst}) and applying
(\ref{AppliedIntByParts}) to the $(a+bd)$-term in (\ref{PreGDotEst}), one obtains
\begin{eqnarray*} 
   \dot{\lambda} \int_N \psi^2~ dV  & \geq & \nonumber
   - 2 b \cdot k \cdot (\lambda+c) \cdot \epsilon^{-1} 
   \int_{I \times M}  \frac{\dot{\rho}}{\rho} \cdot \widehat{\psi}^2~ dV~  \\
   & &  +~ C' \epsilon^{-2a-3} \int_{I \times M} \psi^2  dV~   -~ C'' \cdot (\lambda+1)
    \int_K \psi^2~ dV. \nonumber
\end{eqnarray*}
for some constants $C'$ and $C''$.  Thus, since $a \leq  -1$, there 
exists a positive constant $c'$ such that
\begin{equation*}
   \dot{\lambda}~ \geq~ -2 b  \cdot k \cdot (\lambda+c') \cdot  \epsilon^{-1}~ -~ C''(\lambda+c'),
\end{equation*}
Division by $\lambda+c'$ gives
\begin{equation*}
    \frac{d}{d \epsilon} \log(\lambda+c')~ \geq~ 
    -2 b  \cdot k \cdot   \epsilon^{-1}~ -~ C'',
\end{equation*}
and hence, by arguing as in the proof of Theorem \ref{APrioriEstimate}, one finds
that $\epsilon^{2kb} \cdot \lambda(\epsilon)$ converges as $\epsilon$
tends to zero.  Therefore the claim follows from Theorem \ref{Bootstrap}.
\end{proof}

\medskip

\begin{lem} \label{MinusBound}
There exists $k<1$ such that
\begin{equation}
  \int_{I \times M}  \frac{\dot{\rho}}{\rho} \cdot \rho^{-2b}
    h(\nabla_h \psi_-, \nabla_h \psi_-)~  dV~ \leq~
  k \cdot \lambda \cdot
 \epsilon^{-1}   
 \int_{I \times M} \psi_-^2~ dV.
\end{equation}
\end{lem}
\medskip

\begin{proof}
If $\mu^*$ is less than the smallest positive eigenvalue of $\Delta_h$, 
then $\psi_{-}=0$ and the claim follows. Otherwise, let $\mu_{-}$ be the 
largest $\Delta_h$-eigenvalue that is less than $\mu^*$.  From the 
definition of $\psi_-$ we have
\begin{equation} \label{PsiMinusBound}
\int_{\{t\} \times M} 
  h(\nabla_h \psi_-, \nabla_h \psi_-)~  dV_h~
 \leq~ \mu_- \int_{\{t\} \times M} 
  \psi_-^2~  dV_h
\end{equation}
Since, by hypothesis, $\frac{\mu_-}{\mu^*}<1$, there exists $k<1$ such that
for sufficiently small $\epsilon$ 
\begin{equation}
 \mu_{-}~ \leq~ k \cdot \rho^{2b}(\epsilon, 0)
 \cdot \lambda(\epsilon).
\end{equation}
Hence, using the (global) hypothesis that 
$\rho^{2b}(\epsilon, 0) 
 \leq \rho^{2b}(\epsilon, t)$ 
for all $t$, we have
\begin{equation} \label{Ultimate}
\rho^{-2b}(\epsilon,t) \cdot \mu_{-}~
  \leq~ k \cdot \lambda(\epsilon).
\end{equation}
Combining this with (\ref{PsiMinusBound}) gives
\begin{equation} \label{Rho2b}
\int_{I \times M} \frac{\dot\rho}{\rho}~
    \rho^{-2b} 
h(\nabla_h \psi_-, \nabla_h \psi_-)~  dV_h~
\leq~ k \cdot \lambda 
   \cdot \int_{I \times M} 
  \frac{\dot\rho}{\rho}~ \psi_-^2~  dV_h.
\end{equation}
Therefore, the claim follows from Proposition \ref{MaxedEpsilon}.
\end{proof}

\medskip

\begin{lem} \label{PlusBound}
There exists $k<1$ such that
\begin{eqnarray} \label{PlusFormula} \nonumber
   \int_{I \times M}  \frac{\dot{\rho}}{\rho} \cdot \rho^{-2b}
   h(\nabla_h \psi_+, \nabla_h \psi_+)~ dV~
   &\leq & \left( k \cdot \lambda + C \epsilon^{-2a-2} \right) \cdot
   \epsilon^{-1}  \int_{I \times M} \psi_+^2~ dV~ \\ &&
 +~ C'\cdot (\lambda+1) \int_K \psi^2~ dV.
\end{eqnarray}
\end{lem}

\medskip

\begin{proof}
By (\ref{Gradient}) and (\ref{LocalMetric}) we have
\begin{equation} \label{hTog}
  \rho^{-2b} \cdot h(\nabla_h \psi_+,\nabla_h \psi_+)~ \leq~
   \rho^{2b} \cdot h(\nabla_g \psi_+,\nabla_g \psi_+)~ \leq~
  g(\nabla_g \psi_+,\nabla_g \psi_+).
\end{equation}
Using Lemma \ref{IntByPartsLemma}, one obtains (\ref{AppliedIntByParts})
with $\psi$ replaced by $\psi_+$. From this and (\ref{hTog})
it follows that the left hand side of (\ref{PlusFormula}) is bounded above by
\begin{equation} \label{LeftHandBound}
 \lambda \int_{I \times M} \frac{\dot{\rho}}{\rho} \cdot \psi_+^2~  dV 
 +~ C \cdot \epsilon^{-2a-3} \int_{I\times M} \psi_{+}^2~ +~ C'(\lambda+1)
    \int_{(2I \setminus I) \times M}  \psi_+^2~  dV.
\end{equation} 
Thus it will suffice to show that
\begin{equation} \label{MuPlus}
 \int_{I \times M} 
 \frac{\dot{\rho}}{\rho}~ \psi_+^2~ dV~
 \leq~ k \cdot  \epsilon^{-1}
 \int_{I \times M} 
 \psi_+^2~ dV
\end{equation}
for some $k<1$.

Towards verifying (\ref{MuPlus}), we will apply Lemma \ref{BootstrapLemma}.
Namely, let $\mu_+$ to be the smallest of all eigenvalues that are greater 
than $\mu^*$ and let $\delta= (\mu_{+} -\mu^*)/2$. 
Note that from the definition of $\psi_+$, for each $t \in I$
\begin{equation*} 
 \int_{\{t\}\times  M} 
   \psi_+ \cdot \Delta_h \psi_+~ dV_h~
  \geq~ \mu_+  \int_{\{t\}\times  M}
      \psi_+^2~ dV_h.
\end{equation*}
Note also that (\ref{Dies})  follows in this case from (\ref{mDefn}).
Therefore Lemma \ref{BootstrapLemma} provides $C>0$ such that 
\begin{equation} \label{The52Estimate2}
  \int_{\frac{1}{2} A(\epsilon) \times M} 
  \rho^{-1} \cdot \psi_+^2~ dV_g~
  \leq~  C  \cdot \epsilon^{2b-1} 
  \int_{A(\epsilon) \times M} 
        \psi_+^2~ dV_g.
\end{equation}

By Lemma \ref{DiameterBound}, there exists $\eta>0$ such that
$[-\eta \epsilon , \eta \epsilon ] \subset A(\epsilon)$.
Define
\begin{equation*}
   k'= \sup_{ 2 |s|< \eta } 
   \frac{\dot\rho(1,s)}{ \rho(1,s)}
\end{equation*}
It follows from Proposition \ref{MaxedEpsilon}
that $k' < 1$. By homogeneity, we have
\begin{equation} \label{KPrime}
   \frac{\dot\rho(\epsilon,t)}{\rho(\epsilon,t)}
   \leq~ k' \cdot \epsilon^{-1}
\end{equation}
for all $|t| < 2^{-1} \eta \cdot \epsilon$.
In particular, estimate (\ref{KPrime}) holds for 
all $t \in B(\epsilon)$. Thus, by integrating 
this estimate over $B(\epsilon)$ and combining
with (\ref{The52Estimate2}) we obtain
\begin{equation}
 \int_{I \times M}  \frac{\dot{\rho}}{\rho}~ \psi_+^2~ dV~
  \leq~ (k' + C \epsilon^{2b}) \cdot  \epsilon^{-1} \int_{I \times M} \psi_+^2~ dV
\end{equation}
Therefore, since $b>0$, the claimed (\ref{MuPlus}) is proven.
\end{proof}

\medskip

\medskip



\section{A vacuous case}

\label{SectionMainResult}

\medskip


\begin{thm}
Each eigenvalue branch converges to a finite limit.
\end{thm}
\medskip

\begin{proof}
We assume that the hypothesis of Theorem \ref{ScarTheorem} is not satisfied and derive a 
contradiction. Namely, we assume that $\mu^*$ is a positive $\Delta_h$-eigenvalue and obtain
a contradiction in the form of two conflicting estimates. In particular, it is enough to show that 
there exist constants $C_1, C_2>0$ such that for small $\epsilon$
\begin{equation}  \label{Conflict}
  C_1 \cdot \epsilon^{-b-a-1}~ \leq~
  \lambda(\epsilon)
  -\mu^* \cdot \rho^{-2b}(\epsilon, 0)~ \leq~
  C_2 \cdot \epsilon^{-2a-2}.
\end{equation}
This is impossible since $b>0$ and  $-a-1 \leq 0$ and hence $-b-a-1 < -2a-2$.  
The respective sides of (\ref{Conflict}) are given below 
as Lemma \ref{LeftConflict} and Lemma \ref{RightConflict}.
\end{proof}

\medskip



\begin{lem}[Left Hand Estimate]
\label{LeftConflict}
Suppose that $\mu^*$ is a positive
$\Delta_h$-eigenvalue.
Then there exists a constant $C>0$ such that
for all small $\epsilon$
\begin{equation} \label{ReLeft}
 C \cdot \epsilon^{-b-a-1}~
 \leq~ \lambda   -\mu^* \cdot 
 \rho^{-2b}(\epsilon,0). 
\end{equation}
\end{lem}

\medskip

\begin{proof}
We first claim that it suffices to show that there exists $\delta>0$ such that
for small $\epsilon$
\begin{equation} \label{bMinus1}
  \dot{\lambda}~ \geq~ -2b \cdot  \lambda 
  \cdot \epsilon^{-1}  \cdot \left( 1
  - \delta  \cdot \epsilon^{-a-1} 
  \lambda^{-\frac{1}{2}}  \right).
\end{equation}
Indeed, since $\rho^{2b}(\epsilon, 0) 
   \cdot  \lambda(\epsilon) \rightarrow \mu^*$,
there exists $\delta'>0$ such that
$\lambda^{\frac{1}{2}} \geq
 \delta' \epsilon^{-b}$ for
small $\epsilon$. Hence we would
have
\begin{equation} \label{bMinus12}
  \dot{\lambda}~ \geq~ -2b \lambda 
  \cdot \epsilon^{-1}  
  - \delta \cdot \delta'  
  \cdot \epsilon^{-b-a-2}.
\end{equation}
Note that since
$(\dot{\rho}/\rho)(\epsilon, 0)
       = \epsilon^{-1}$
\begin{equation}
\frac{d}{d \epsilon} \left( 
 \rho^{2b}(\epsilon, 0) 
        \cdot \lambda(\epsilon) \right)~
  \geq~ \rho^{2b}(\epsilon,0) \cdot 
  \left(2b \cdot \lambda \cdot \epsilon^{-1} 
  + \lambda'(\epsilon) \right).
\end{equation}
Thus, since $\rho^{2b}(\epsilon,0)
= \rho(1,0)\cdot  \epsilon^{2b}$ 
it would 
follow from (\ref{bMinus12}) that there
exists a constant $\delta''>0$ such that
   \begin{equation} \label{Product}
   \frac{d}{d \epsilon} \left( 
   \rho^{2b}(\epsilon, 0) 
   \cdot \lambda(\epsilon) \right)~
   \geq~ \delta'' \cdot \epsilon^{b-a-2}
\end{equation}
Note that since $a \leq -1$ and $b>0$,
we have $b-a-2>-1$. Thus, since 
$\lim_{\epsilon \rightarrow 0} 
 \rho^{2b}(\epsilon, 0) \cdot
  \lambda(\epsilon)=\mu^*$, we could then
integrate  (\ref{Product}) over
$[0, \epsilon]$ and would find that 
\begin{equation} \label{PreLeft}
\rho^{2b} (\epsilon,0) 
\cdot \lambda(\epsilon)~
- \mu^*~ \geq~  \delta'' \cdot \epsilon^{b-a-1}.
\end{equation}
Since $\rho^{2b}(\epsilon,0)=
 \rho^{2b}(1,0) \cdot \epsilon^{2b}$,
we would then obtain (\ref{ReLeft}) by 
dividing both sides of 
(\ref{PreLeft}) 
by $\rho^{2b}(\epsilon,0)$.

Recall from (\ref{PsiStar}) that $\psi^*$ denotes the
fibrewise projection of $\psi$ onto  the $\mu^*$ eigenspace of $\Delta_h$. 
Letting $\bar{\psi} = \psi_+ + \psi_-$, we have 
$\hat{\psi}=  \psi^* +  \overline{\psi}$.
We claim that to verify (\ref{bMinus1}) it suffices to show that 
\begin{equation} \label{EpsilonInterval}
  \int_{ (\epsilon I) \times M}  \frac{ \dot{\rho} }{\rho}
  \cdot (\psi^*)^2~ dV~  \leq~  r(\epsilon) \cdot \epsilon^{-1}
 \cdot \int_{(\epsilon I) \times M} (\psi^*)^2~ dV
\end{equation}
where  
\begin{equation} \label{RDefn}
r(\epsilon)=1- \delta'  \cdot \epsilon^{-a-1} \lambda^{-\frac{1}{2}}
\end{equation}
and $\delta'>0$.
To see this, note that by Proposition \ref{MaxedEpsilon}, there exists $k<1$ such that 
$(\dot{\rho}/\rho)(\epsilon,t) \leq k \cdot \epsilon^{-1}$ for $|t| \geq \epsilon$. 
It follows that
\begin{equation} \label{FullInterval}
  \int_{(I\setminus \epsilon I) \times M} \frac{\dot{\rho}}{\rho} \cdot (\psi^*)^2~ dV~
  \leq~  k \cdot \epsilon^{-1} \cdot \int_{(I\setminus \epsilon I) \times M}
  (\psi^*)^2~ dV.
\end{equation} 
Since $\mu^*>0$, we have $\lambda \rightarrow \infty$, and hence 
since $a+1 \leq 0$, for any $k<1$, we have $r(\epsilon)>k$
for all sufficiently small $\epsilon$.  Thus, by combining 
(\ref{EpsilonInterval}) and (\ref{FullInterval}), we would obtain
 \begin{equation} \label{FullerInterval}
  \int_{I \times M} \frac{\dot{\rho}}{\rho} \cdot (\psi^*)^2~ dV~
  \leq~  r(\epsilon) \cdot \epsilon^{-1} \cdot \int_{I \times M}
  (\psi^*)^2~ dV
\end{equation} 
for $\epsilon$ small. Applying the argument in (\ref{hTog}) and 
(\ref{LeftHandBound}) with $\psi^+$ replaced by $\psi^*$, we would have 
\begin{eqnarray} \label{StarInterval} \nonumber
  \int_{I \times M} \frac{\dot{\rho}}{\rho} \cdot \rho^{-2b} \cdot 
  h(\nabla_h \psi^*,\nabla_h \psi^*) dV~
  &\leq& \left( \lambda \cdot  r(\epsilon) + C \epsilon^{-2a-2}\right)
  \cdot \epsilon^{-1} \cdot \int_{I \times M} (\psi^*)^2~ dV \\ 
  && +~ C' \cdot (\lambda +1) \int_K \psi^2 dV
\end{eqnarray} 
(In this and what follows, $C$ and $C'$ represent generic constants.)
Recall that in (\ref{AppliedhParseval}) we had $\widehat{\psi}= \psi_+ + \psi_-$,
and hence we have (\ref{AppliedhParseval}) with $\widehat{\psi}$ replaced by $\overline{\psi}$:
 \begin{eqnarray} \label{NonStar} \nonumber
  \int_{I \times M} \frac{\dot{\rho}}{\rho} \cdot 
   h(\nabla_h\overline\psi,\nabla_h \overline\psi) dV~
  &\leq&  (\lambda \cdot k +C \cdot \epsilon^{-2a-2}) 
  \cdot \epsilon^{-1} \cdot \int_{I \times M}
  (\overline\psi)^2~ dV~  \\ && +~  C' \cdot (\lambda+1) \int_K \psi^2 dV.
\end{eqnarray}
As pointed out above, $r(\epsilon)>k$ for small $\epsilon$.
Hence by applying Parseval's principle as in (\ref{hParseval}) 
and (\ref{PlainParseval}) to $\widehat{\psi}= \psi^* +\overline{\psi}$,
we could combine (\ref{StarInterval}) and (\ref{NonStar}) to find that 
\begin{eqnarray} \label{AppliedhParseval2} \nonumber
    \int_{I \times M} \frac{\dot{\rho}}{\rho} \cdot \rho^{-2b} 
    \cdot  h(\nabla_h \widehat \psi, \nabla_h  \widehat \psi)~ dV~
    &\leq&  \left( \lambda \cdot r(\epsilon)  +C \cdot \epsilon^{-2a-2} \right)
    \cdot \epsilon^{-1} \int_{I \times M} \widehat\psi^2~ dV~ \\ &&
    +~    C' \cdot (\lambda+1) \int_K \psi^2~ dV.
\end{eqnarray}
Combining this with (\ref{bPosBound})  would yield $k<1$ such that
\begin{equation} \label{PsiBarEstimate2}
 \int_{I \times M} \dot{g}(\nabla \widehat{\psi}, \nabla \widehat{\psi})~ dV~ \leq~
 \left( 2b \cdot \lambda \cdot r(\epsilon) +C \cdot \epsilon^{-2a-2}  \right) \cdot  \epsilon^{-1}   
 \int_{I \times M} \widehat\psi^2~ dV~ +~ C'\cdot (\lambda+1) \int_K \psi^2 dV.
\end{equation}
where $C$ denotes a (generic) positive constant.
By substituting (\ref{PsiBarEstimate2}) into (\ref{PreGDotEst}) and applying
(\ref{AppliedIntByParts}) to the $(a+bd)$-term in (\ref{PreGDotEst})
and using (\ref{HatReduction}), one would obtain
\begin{equation} \label{WithC}
   \dot{\lambda} \int_N \psi^2~ dV~ \geq~ - \left(2 b \cdot \lambda \cdot 
   r(\epsilon)+ C \epsilon^{-2a-2} \right) \cdot  \epsilon^{-1} 
   \int_{I \times M} \psi^2 dV~  -~ C' (\lambda+1) \int_{K} \psi^2~ dV.
\end{equation}
Note that $C'(\lambda+1) < 2b \lambda \cdot r \cdot \epsilon^{-1}$
for small $\epsilon$. Hence  from (\ref{WithC}) one would have
\begin{equation} \label{WithC2}
\dot{\lambda} \int_N \psi^2~ dV~ \geq~ - \left(2 b \cdot \lambda \cdot 
   r(\epsilon)+ C \epsilon^{-2a-2} \right) \cdot  \epsilon^{-1} 
   \int_{N} \psi^2 dV~ 
\end{equation}
for small $\epsilon$.  Since $\lambda \rightarrow \infty$ and $-a-1 \geq 0$
\begin{equation} 
 C \cdot \epsilon^{-2a-2}~ \leq~  \frac{\delta'}{2} \cdot \epsilon^{-a-1} 
\end{equation}
for small $\epsilon$. Thus by choosing $\delta= \delta'/2$ and recalling (\ref{RDefn})
we would obtain (\ref{bMinus1}) from (\ref{WithC2}). And hence (\ref{bMinus1}) follows
from (\ref{EpsilonInterval}).

As a first step toward the verification of (\ref{EpsilonInterval}), we  
rescale in $\epsilon$. In particular, let $h(s)= \dot{\rho} \rho^{-1}(1,s)$,
and $\sigma(s)= \rho^{2a}(1,s)$, and for each $\epsilon>0$ and $m \in M$, define
\begin{equation} \label{VDefn}
 v_{\epsilon,m}(s)=
 \rho^{\frac{-a+bd}{2}}(1,s) 
\cdot \psi^*(\epsilon s, m).
\end{equation}
Then using homogeneity, we find that
\begin{equation} \label{Reduction1}
 \int_{(\epsilon I) \times M}  \frac{\dot{\rho}}{\rho} \cdot (\psi^*)^2~  dV_g~
 =~ \epsilon^{a+bd} \int_{M} \int_{I} h(s) \cdot v_{\epsilon,m}(s)^2 \cdot  \sigma(s)~
    ds~ dV_h
\end{equation}
and
\begin{equation}  \label{Reduction2}
 \int_{(\epsilon I) \times M}
 (\psi^*)^2~ 
 dV_g~
 =~
 \epsilon^{a+bd+1} \int_{M} \int_{I}
v_{\epsilon,m}(s)^2
 \cdot  \sigma(s)~
 ds~ dV_h.
\end{equation}
Since $\Delta \psi^* = \lambda \psi^*$
and since
$\psi^*(t,\cdot)$ belongs to the
$\mu^*$ eigenspace of $\Delta_h$,
we have $-L \psi^*+ \mu^* \rho^{-2b}
\psi^* = \lambda \psi^*$ from 
(\ref{Laplacian}).
It follows from (\ref{VDefn}) and 
homogeneity
that $v_{\epsilon,m}$ satifies the 
following ordinary differential equation
\begin{equation} \label{VODE}
   v''(s)~
=~ \eta \cdot \rho^{2a-2b}(1,s) \cdot \left( \frac{\mu^*}{\lambda \cdot \epsilon^{2b}}
    -\rho^{2b}(1,s) \right) \cdot v(s) + g(s) \cdot v(s)
\end{equation}
where $\eta= \epsilon^{2a+2} \cdot \lambda(\epsilon)$
and $g(s)$ is a bounded smooth function.  Hence, by (\ref{Reduction1})
and (\ref{Reduction2}), to prove (\ref{EpsilonInterval}) it 
will suffice to prove that there exists $\delta>0$ such that for any 
solution $v$ to (\ref{VODE}) we have
\begin{equation} \label{DesireV}
  \int_{I}  h(s) \cdot v(s)^2  \cdot  \sigma(s)~ ds~ 
  \leq~ \left(1- \delta \cdot  \eta^{-\frac{1}{2}} \right)
  \int_{I}  v(s)^2  \cdot  \sigma(s)~ ds.
\end{equation}

Towards verification of (\ref{DesireV}) we apply Lemma \ref{ScarLemma} to (\ref{VODE}). 
Indeed, by hypothesis $\mu^* \cdot \lambda^{-1} \cdot \epsilon^{-2b}$ is
bounded, and hence Lemma \ref{ScarLemma} applies to give a constant $C>0$
such that 
\begin{equation}
 \int_{ \eta^{-\frac{1}{4}}I } v(s)^2 \cdot \sigma(s)~ ds~
 \leq~ C \int_{ \eta^{-\frac{1}{4}}(2I \setminus I)}  v(s)^2 \cdot \sigma(s)~ ds.
\end{equation}
It follows that 
\begin{equation} \label{PrimeNoPrime}
 \int_{I} v(s)^2 \cdot \sigma(s)~ ds~ 
 \leq~ (C+1) \int_{I \setminus (2 \eta^{-\frac{1}{4}}I) } v(s)^2 \cdot \sigma(s)~ ds.
\end{equation}

By Proposition \ref{MaxedEpsilon} and the strict convexity of $\rho$, 
there exists $\delta'>0$ such that 
$h(s)= \dot{\rho} \rho^{-1}(1,s) \leq 1- \delta' \cdot s^2$ for $s \in I$. 
It follows that
\begin{equation*}
 \int_{I} h(s) \cdot v(s)^2 \cdot  \sigma(s)~ ds~ \leq~
 \int_{I} v(s)^2 \cdot  \sigma(s)~ ds~ 
 -~ 4 \delta' \cdot \eta^{-\frac{1}{2}} 
 \int_{I \setminus (2 \eta^{-\frac{1}{4}}I)  } v(s)^2 \cdot  \sigma(s)~ ds.
\end{equation*}
By substituting (\ref{PrimeNoPrime}) and choosing $\delta= 4 \delta' \cdot (1+C)^{-1}$
we obtain (\ref{DesireV}). The proof is complete.
\end{proof}

\medskip



\begin{lem}[Right Hand Estimate] 
\label{RightConflict}
Suppose that $\mu^*$ is a positive $\Delta_h$-eigenvalue.
Then there exists a constant $C$ such that for small $\epsilon>0$
\begin{equation} \label{RightEst}
  f(\epsilon)~ =~  \lambda(\epsilon)~ -~
  \mu^* \cdot \rho^{-2b}(\epsilon,0)~  
  \leq~ C \epsilon^{-2a-2}.
\end{equation}
\end{lem}
\medskip

\begin{proof}
It will suffice to show that
\begin{equation}  \label{LogA}
   \frac{d}{d \epsilon} 
   f(\epsilon)~ \geq~ 
  -(2a \cdot f(\epsilon)+C \epsilon^{-2a-2}) 
   \cdot \epsilon^{-1}.
\end{equation}
Indeed, we may suppose that $2 \cdot f(\epsilon) \geq C \epsilon^{-2a-2}$, 
for otherwise we are done.  Hence (\ref{LogA}) implies
\begin{equation*}  
   \frac{d}{d \epsilon} f(\epsilon)~ \geq~ 
  -(2a+2) \cdot f(\epsilon)   \cdot \epsilon^{-1}.
\end{equation*}
By Lemma \ref{LeftConflict}, we have $f>0$ for small $\epsilon$,  
and thus division would give
\begin{equation*}
   \frac{d}{d \epsilon} \log \left(  \epsilon^{2a+2}  \cdot f(\epsilon) \right)~ 
   \geq~   0.
\end{equation*}
By integrating over $[\epsilon, \delta]$
with $\delta$ small we would find that
\begin{equation*}
 \log \left( \frac{\delta^{2a+2} \cdot f(\delta)}{ \epsilon^{2a+2} \cdot f(\epsilon)}  \right)
  \geq 0.
\end{equation*}
Exponentiation would then give (\ref{RightEst}) with $C= \delta^{2a+2} \cdot f(\delta)$.

To verify (\ref{LogA}) we will estimate $\lambda'(\epsilon)$ using
the methods of \S \ref{SectionLowerBound}. From (\ref{GDotFormula}) 
and (\ref{Gradient}) we have
\begin{equation} \label{RhoDotApplied}
  \dot{g}(\nabla \widehat \psi, \nabla \widehat \psi)~
  =~ 2a \cdot \frac{\dot{\rho}}{\rho} \cdot  \rho^{-2a} \cdot (\partial_t \widehat \psi)^2~
  +~ 2b \cdot  \frac{\dot{\rho}}{\rho} \cdot \rho^{-2b} \cdot h(\nabla_h \widehat \psi, \nabla_h \widehat \psi).
\end{equation}
By substituting (\ref{GDecomposed}) into (\ref{AppliedIntByParts})
one obtains that
\begin{eqnarray} \label{PartsAppliedTwo} \nonumber
   \left|  \int_{I \times M}  \frac{\dot{\rho}}{\rho} \cdot \rho^{-2a} (\partial_t \widehat \psi)^2 dV
  + \int_{I \times M}  \frac{\dot{\rho}}{\rho} \cdot \rho^{-2b} \cdot 
                     h(\nabla_h \widehat \psi, \nabla_h \widehat \psi) dV    
  - \lambda \int_{I \times M}  \frac{\dot{\rho}}{\rho}  \widehat \psi^2 dV \right|~  &&  \\
  \ \ \ \ \leq  C \int_{I\times M} \rho^{-2a-3} \cdot \widehat \psi^2~ dV~ 
          +  C'(\lambda + 1) \int_{K}  \widehat \psi^2~ dV  &&
\end{eqnarray}
Thus, by integrating (\ref{RhoDotApplied}) and using (\ref{PartsAppliedTwo})
we find that 
\begin{eqnarray} \label{Begin} \nonumber
  \int_{I\times M} \dot{g}(\nabla \widehat \psi, \nabla \widehat \psi) dV
  &\leq&   2 a \cdot \lambda \int_{I\times M} \frac{\dot{\rho}}{\rho}~ \widehat \psi^2~dV  \\
  && + (2b-2a) \cdot \int_{I\times M} \frac{\dot{\rho}}{\rho}~ \rho^{-2b} \cdot
  h(\nabla_h \widehat \psi, \nabla_h \widehat \psi)~dV  \\
  & &+  C \cdot \epsilon^{-2a-3} \int_{I\times M} \widehat\psi^2~ dV   
  + C'(\lambda +1)  \int_{(2I\setminus I) \times M}  \widehat \psi^2~ dV. \nonumber
\end{eqnarray}
We may estimate the righthand side of (\ref{Begin}) by splitting the integral over the sum
$\widehat{\psi}= \psi^* +\overline{\psi}$ where $\overline{\psi}=\psi_+ + \psi_-$. 
To be precise,  apply Parseval's principle---as in (\ref{hParseval})---to find that
\begin{equation*} 
  \int_{\{t\} \times M} h(\nabla_h \widehat \psi, \nabla_h \widehat  \psi)~ dV_h~
  =~  \int_{\{t\} \times M} h(\nabla_h \psi^*, \nabla_h \psi^*)~ dV_h~
  +~  \int_{\{t\} \times M} h(\nabla_h \overline{\psi}, \nabla_h \overline{\psi})~ dV_h
\end{equation*}
as well as 
\begin{equation} \label{BabyParseval}
    \int_{\{t\} \times M} (\widehat \psi)^2~dV_h~ =~ \int_{\{t\} \times M} (\psi^*)^2~dV_h~ 
   +~ \int_{\{t\} \times M} (\overline \psi)^2~ dV_h.
\end{equation}
It follows that if we let $E(\widehat{\psi})$ denote the 
right hand side of (\ref{Begin}), then
\begin{equation} \label{ESeparate}
  E(\widehat{\psi})~ =~  E(\psi^*)~ +~ E(\overline{\psi}).
\end{equation}

By definition $\Delta_h \psi^* = \mu^* \psi^*$ and hence  
$\int_{M} h(\nabla_h \psi^*, \nabla_h \psi^*) dV_h= \mu^* \int_M (\psi^*)^2 dV_h$.
It follows that
\begin{equation}  \label{EStar}
 E(\psi^*)~ \leq~    \left(q(\epsilon)+ C \epsilon^{-2a-2} \right)
         \cdot \epsilon^{-1}  \int_{I\times M} (\psi^*)^2~ dV  
      + C'(\lambda +1)  \int_{K}  \psi^2~ dV.
\end{equation}
where
\begin{equation} \label{QDefn}
 q(\epsilon)=2 a \cdot \lambda(\epsilon) +
 (2b-2a)\cdot \mu^* \cdot \rho^{-2b}(\epsilon,0).
\end{equation}

To estimate $E(\overline{\psi})$, apply (\ref{PartsAppliedTwo})  `in reverse' 
to find that 
\begin{eqnarray*}
 E(\overline{\psi})
  &\leq& ~ 2a  \int_{I \times M} 
  \frac{\dot{\rho}}{\rho} \cdot  \rho^{-2a} \cdot (\partial_t \overline \psi)^2~ dV~
  +~ 2b    \int_{I \times M} 
 \frac{\dot{\rho}}{\rho} \cdot \rho^{-2b} 
   \cdot h(\nabla_h \overline \psi, \nabla_h \overline \psi)~ dV   \\ 
  & &+  C'' \cdot \epsilon^{-2a-3} \int_{I\times M} \overline\psi^2~ dV   
  + C'''(\lambda +1)  \int_{(2I\setminus I) \times M}  \overline \psi^2~ dV. \nonumber
\end{eqnarray*}
Thus by using the assumption $a<0$, the formula (\ref{hParseval}),
and  Lemmas \ref{PlusBound} and \ref{MinusBound}, we obtain $k<1$ such that 
\begin{equation} \label{EOverline}
   E(\overline{\psi})~  \leq~  (2k \cdot b \lambda +  C'' \cdot \epsilon^{-2a-2})
   \cdot  \epsilon^{-1} \int_{I \times M}  \overline{\psi}^2~ dV
  + C'''(\lambda +1)  \int_{K}  \overline \psi^2~ dV
\end{equation}

By combining (\ref{HatReduction}), 
(\ref{Begin}), (\ref{ESeparate}), (\ref{EStar}), and (\ref{EOverline}),
we find (generic) constants $C, C'$ such that
\begin{eqnarray*}
\int_{I \times M} \dot{g}(\nabla \psi, \nabla \psi)~ dV~
&\leq& (q + C \epsilon^{-2a-2})\epsilon^{-1}\int_{I \times M} (\psi^*)^2~ dV~\\
&&  +~  (2kb \lambda + C \epsilon^{-2a-2})\epsilon^{-1}\int_{I \times M} (\overline{\psi})^2~ dV~
  +~ C'(\lambda +1)  \int_{K}  \psi^2~ dV.
\end{eqnarray*}
Hence by substituting this into (\ref{PreGDotEst}) and applying (\ref{AppliedIntByParts})
to the middle term in (\ref{PreGDotEst}), we obtain 
\begin{eqnarray*}
\dot{\lambda} \int_{N} \psi^2~ dV~
&\geq& -(q + C \epsilon^{-2a-2}) \epsilon^{-1}\int_{I \times M} (\psi^*)^2~ dV~\\
&&  -~  (2kb \lambda + C \epsilon^{-2a-2})\epsilon^{-1}
   \int_{I \times M} (\overline{\psi})^2~ dV~
  -~ C'(\lambda +1)  \int_{K}  \psi^2~ dV
\end{eqnarray*}
for possibly different (generic) constants. By using
(\ref{BabyParseval}), we find that
\begin{eqnarray} \label{StarEstimate} \nonumber
\dot{\lambda} \int_{N} \psi^2~ dV~ 
  &\geq& -\left(\max\{q,2kb \lambda\}  + C \epsilon^{-2a-2}\right)
\epsilon^{-1}  \int_{I \times M} \widehat \psi^2~ dV~ \\
 && -~ C'(\lambda+1)  \int_{K}  \psi^2~ dV.
\end{eqnarray}

We claim that there exists $\epsilon_0>0$ such that for all
$\epsilon \leq \epsilon_0$ 
\begin{equation} \label{QBig}
    q(\epsilon)~ \geq~ 2 b \cdot k \cdot \lambda(\epsilon).
\end{equation}
for all small $\epsilon$.  To see this, first note that we cannot have
$q <  2 b \cdot k \cdot \lambda$ for all small $\epsilon$. For then by 
(\ref{StarEstimate}) we would have that
$\dot{\lambda} \geq -(2b k \lambda + C)\cdot \epsilon^{-1}  -C' \lambda$
for all small $\epsilon$. Using the argument of Theorem \ref{APrioriEstimate}, we would
then find that $\epsilon^{2kb} \lambda$ is bounded and hence $\lambda$ 
would converge by Theorem \ref{Bootstrap}. This would contradict the assumption
that $\mu^*>0$.

Hence there exists $\epsilon_0>0$ such that (\ref{QBig}) is true for
$\epsilon= \epsilon_0$. A calculation shows that (\ref{QBig})  holds if and only if
\begin{equation} \label{Under}
 \left(  \frac{1-k}{-a+b}  \right)~
\geq~ \frac{f(\epsilon)}{\lambda(\epsilon) }.
\end{equation}
Thus to prove (\ref{QBig}) for $\epsilon \leq \epsilon_0$, it 
suffices to show that $f/\lambda$ is increasing for small $\epsilon$. 
To see that this is true, note that by homogeneity 
$\partial \rho^{-2b}(\epsilon,0) = -2b \cdot
  \rho^{-2b}(\epsilon,0) \cdot \epsilon^{-1}$,
and hence
\begin{equation} \label{FPrime}
  f'(\epsilon) = \lambda'(\epsilon) 
  + 2b\mu^* \cdot \rho^{-2b}(\epsilon,0) 
  \cdot \epsilon^{-1}.
\end{equation}
By Lemma \ref{LeftConflict},  $f(\epsilon)>0$
for small $\epsilon$. Thus from (\ref{FPrime})
we have $\partial \log(f(\epsilon)) \geq \partial
    \log(\lambda(\epsilon))$, and hence
$\partial (f/\lambda) \geq 0$ as desired.

Therefore,  since (\ref{QBig}) holds true,  
(\ref{StarEstimate}) yields
\begin{equation}  \label{FirstDotB2}
   \dot{\lambda} \int_N \psi^2~ dV~
   \geq~ -\left(q(\epsilon) 
   + C \cdot \epsilon^{-2a-2} \right)
   \cdot \epsilon^{-1}
   \int_{I \times M} \widehat{\psi}^2~ dV~
  -~    C' (\lambda+1)
  \int_{K} 
  \psi^2~ dV.
\end{equation}
Moreover, for small $\epsilon$ we have $2 k' b  \cdot \epsilon^{-1} \geq 2 C'$,
and hence by (\ref{StarEstimate}) we have $q > C' (\lambda+1)$. 
it follows from 
(\ref{FirstDotB2}) that 
\begin{equation*}  
   \dot{\lambda} 
   \geq~ -\left(q(\epsilon) 
   + C \cdot \epsilon^{-2a-2} \right)
   \cdot \epsilon^{-1}.
\end{equation*}
Into this substitute (\ref{QDefn}), and use both (\ref{FPrime})
and the definition of $f$ in (\ref{RightEst}) to obtain (\ref{LogA}) as desired.
\end{proof}

\medskip


\appendix

\section{On the width of scars}

\label{AppendixScar}
\medskip

Let $f>0$, $g$ and $h$, be continuous 
functions on ${\mathbb R}$. Let $I$
be an interval containing $0$, and let
$\beta \in {\mathbb R}$ and
$\eta \in {\mathbb R}^+$.
Consider the ordinary differntial 
equation
\begin{equation} \label{AbstractODE}
 w''(s) =  \eta \cdot f(s)
 \cdot \left( \beta + s^2 \cdot h(s)
 \right)
 \cdot w(s) 
  + g(s) \cdot w(s).
\end{equation}

\medskip

\begin{lem} \label{ScarLemma}
Let $\sigma>0$ 
be a positive continuous function 
on ${\mathbb R}$.
There exists a constant $C>0$ such 
that for any solution $w$ to 
(\ref{AbstractODE})
\begin{equation} \label{NormComparison}
 \int_{\lambda^{-\frac{1}{4}} I} 
 w^2(s)~ \sigma(s)~ds~
 \leq~ C \int_{\lambda^{-\frac{1}{4}}
 (2 I \setminus I)}   w^2(s)~ \sigma(s)~
 ds. 
\end{equation}
\end{lem} 
\begin{proof}
Let $x(u)=w(\eta^{-\frac{1}{4}}u)$.
Then (\ref{NormComparison}) is equivalent
to
\begin{equation} \label{NormComparison2}
 \int_{I}
 x^2(u)~ \sigma(\eta^{-\frac{1}{4}}u)~du~
 \leq~ C \int_{
 2 I \setminus I}   x^2(u)~ 
 \sigma(\eta^{-\frac{1}{4}}u)~
 du. 
\end{equation}
From (\ref{AbstractODE}), we see that
$x$ satisfies
\begin{equation} \label{AbstractODE2}
 x''(u) =  \eta^{\frac{1}{2}}\cdot
 \beta \cdot f(\eta^{-\frac{1}{4}}u)
 \cdot  x(u)
 + \left( 
 u^2 \cdot h(\eta^{-\frac{1}{4}}u)
 + g(\eta^{-\frac{1}{4}}u) \right)
 \cdot x(u).
\end{equation}
If $\eta^{\frac{1}{2}} \cdot \beta$
is large and positive, then $x(u)$
is uniformly convex on $I$, and 
(\ref{NormComparison2}) follows.
If $\eta^{\frac{1}{2}} \cdot \beta$
is large and negative, then $x(u)$
oscillates rapidly. In particular, 
by \cite{CrnHlb} Chapter V \S 4, 
$x$ differs
from a solution to $y''(u)+ 
(\eta^{\frac{1}{2}}  \beta) \cdot  y(u)$
in $C^0$-norm by order 
$(\eta^{\frac{1}{2}} 
        \beta)^{-\frac{1}{2}}$.
A straightforward calculation shows that
(\ref{NormComparison2}) holds uniformly
for $y$ and hence $w$ for
$\eta^{\frac{1}{2}} \beta$
positive and sufficiently large.
For $\eta^{\frac{1}{2}} \beta$ bounded, 
the claim follows from 
continuous dependence on parameters,
the linearity of the equation, 
and the fact that the $L^2$-norm 
cannot vanish on any nontrivial interval.
\end{proof}


\end{document}